\documentclass[11pt, a4paper]{article}
\usepackage[utf8]{inputenc}
\usepackage{tikz}
\usepackage[T1]{fontenc}
\usepackage[english]{babel}
\usepackage{bbm}
\usepackage{array}
\newcolumntype{M}[1]{>{\centering\arraybackslash}m{#1}}
\usepackage{listings}
\usepackage{ragged2e}
\usepackage{graphicx} 

\usepackage{soul} 

\newcommand{\pcolon}{\mathbin{:}}

\usepackage{breqn}

\usepackage{tikzpagenodes}

\usepackage[font=footnotesize,skip=0pt]{caption}

\usepackage{pgfplots}
\pgfplotsset{width=7cm,compat=1.9}

\usepackage{xcolor}
\colorlet{myorange}{orange!80!black} 
\colorlet{myred}{red!80!black}
\colorlet{mylessred}{red!60!black}
\colorlet{myblue}{blue!80!black}
\colorlet{mygreen}{green!60!black}
\colorlet{mydarkgreen}{green!40!black}
\colorlet{mydarkgray}{gray!70!black}
\colorlet{mygrayHeader}{gray!40!white}

\usepackage{datetime}

\usepackage{lipsum}       
\usepackage{changepage}   

\usepackage[toc,header, page]{appendix}

\usepackage[,linesnumbered,boxed,lined,commentsnumbered]{algorithm2e}

\usepackage{multirow}
\usepackage{float}

\usepackage{hyperref}
\definecolor{Black}{rgb}{0,0,0}
\definecolor{CiteColor}{rgb}{0,0,0}
\hypersetup{
	colorlinks   = true, 
	urlcolor     = black, 
	linkcolor    = myblue, 
	citecolor   = black 
}

\usepackage{natbib}

\usepackage{textcomp}
\usepackage{a4wide}
\usepackage{dsfont}
\usepackage{amsmath,amsthm,amsfonts,amssymb,mathrsfs,graphicx}

\usepackage{enumerate}
\usepackage{enumitem}
\usepackage{aligned-overset}

\usepackage{makeidx}

\usepackage[left=2.5cm,right=2.5cm,top=2.5cm,bottom=2.5cm]{geometry}

\usepackage{cleveref} 

\theoremstyle{plain}
\newtheorem{thm}{Theorem}[section] 

\theoremstyle{definition}

\newtheorem{remx}[thm]{Remark} 
\newenvironment{rem}
{\pushQED{\qed}\remx}
{\popQED\endremx}

\newtheorem{lemma}[thm]{Lemma}

\newtheorem{examplex}[thm]{Example} 
\newenvironment{example}
{\pushQED{\qed}\examplex}
{\popQED\endexamplex}

\newtheorem{assumption}[thm]{Assumption}
\newtheorem{definition}[thm]{Definition}
\newtheorem{proposition}[thm]{Proposition}

\newcommand{\enumthm}{(\alph*)}

\newcommand{\enumproposition}{(\alph*)}

\newcommand{\enumintext}{\arabic*.}

\newcommand{\enumassumption}{(\Roman*)}
\newcommand{\salg}{\mathcal{F}}
\newcommand{\proc}{Z}
\newcommand{\law}{\mu}
\newcommand{\sspa}{\mathcal{Z}}
\newcommand{\dspa}{\mathcal{X}}
\newcommand{\pspa}{\Omega}
\newcommand{\rate}{r}
\newcommand{\ctsset}{\mathcal{C}}
\newcommand{\costset}{\ctsset(\dspa\times\sspa)}
\newcommand{\decset}{\ctsset(\sspa,\dspa)}
\newcommand{\uset}{\ctsset(\sspa)}

\newcommand{\setmap}{\Upsilon}
\newcommand{\off}{u}
\newcommand{\dec}{X}

\newcommand{\decplug}{X^{\mathsf{Plug}\text{-}\mathsf{in}}}
\newcommand{\supfct}{G}
\newcommand{\modsupfct}{H}

\newcommand{\bbE}{\mathbb{E}}
\newcommand{\utheta}{\theta'}
\newcommand{\realinf}{[-\infty,\infty]}
\newcommand{\optimal}[1]{{#1}^\star}

\newcommand{\indfctset}[1]{\mathbbm{1}_{#1}}
\newcommand{\supb}[2]{\sup_{#1}\left\{#2\right\}}

\newcommand{\maxb}[2]{\max_{#1}\left\{#2\right\}}
\newcommand{\phisub}[2]{\phi_{#1}\left[#2\right]}
\newcommand{\phiupsub}[3]{\phi^{\mathsf{#1}}_{#2}\left[#3\right]}
\newcommand{\interior}[1]{\mathring{#1}}

\newcommand{\ratefct}{\mathcal{I}}
\newcommand{\ent}{\ratefct^\mathsf{ENT}}
\newcommand{\rent}{\ratefct^\mathsf{RENT}}
\newcommand{\pseudoinv}[1]{#1^{-1}}
\newcommand{\inv}[1]{#1^{-1}}

\newcommand{\seq}[2]{(#1_{#2})_{#2\in\nat}}
\newcommand{\set}[2]{\left\{#1 \mid #2 \right\}}
\newcommand{\converge}{\to}
\newcommand{\defeq}{\vcentcolon=}

\newcommand{\norm}[1]{\left\| #1 \right\|}
\newcommand{\abs}[1]{\left\lvert#1\right\rvert}
\newcommand{\real}{\mathbb{R}}
\newcommand{\nat}{\mathbb{N}}

\newcommand{\expecsub}[2]{\mathbb{E}_{#1}\left[#2\right]}
\newcommand{\subexpecsub}[2]{\overline{\mathbb{E}}_{#1}\left[#2\right]}

\newcommand{\probsub}[2]{\mathbb{P}_{#1}\left[#2\right]}

\newcommand{\map}{\rightarrow}

\newcommand{\bbP}{\mathbb{P}}

\DeclareMathOperator{\regret}{\mathsf{Regret}}

\DeclareMathOperator{\subjectto}{s.t.}

\DeclareMathOperator*{\argmin}{arg\,min}

\DeclareMathOperator*{\Argmin}{Arg\,min}

\DeclareMathOperator{\dom}{dom}

\newcommand{\defformat}[1]{{\emph{#1}}}

\newcommand{\proofheading}[1]{\textbf{#1}}

\numberwithin{equation}{section}

\title{\Huge Asymptotic Optimality in Data-Driven Decision Making}
\author{Radek Sala\v{c}$^1$, Michael Kupper$^2$, Tobias Sutter$^{3}$}
\date{\small{
    $^1$Department of Computer Science, University of Konstanz, radek.salac@uni-konstanz.de\\
    $^2$Department of Mathematics, University of Konstanz, {kupper@uni-konstanz.de}\\
    $^3$Department of Economics, University of St.Gallen, {tobias.sutter@unisg.ch}\\%
    [2ex]}
  \today
}

\usepackage{glossaries-extra}
\setabbreviationstyle[acronym]{long-short}
\glssetcategoryattribute{acronym}{nohyperfirst}{true}

\newacronym{SAA}{SAA}{sample average approximation}
\newacronym{DRO}{DRO}{distributionally robust optimization}
\newacronym{LDP}{LDP}{large deviation principle}
\newacronym{LP}{LP}{Laplace principle}
\newacronym{usc}{usc}{upper semi-continuous}
\newacronym{lsc}{lsc}{lower semi-continuous}
\newacronym{LLN}{LLN}{law of large numbers}
\newacronym{RLLN}{RLLN}{robust law of large numbers}

\usepackage{acronym}
\renewcommand{\ac}[1]{\gls{#1}}
\renewcommand{\acf}[1]{\gls{#1}}

\begin{document}

\maketitle
\begin{center}
    {\bf Abstract}
\end{center}
Given data generated by an observable stochastic process, we study how to construct statistically optimal decisions for general stochastic optimization problems. Our setting encompasses non-standard data structures, including data originating from heterogeneous sources or from randomly evolving data-generating mechanisms. We propose a decision-making approach that identifies optimal decisions for which a specific notion of risk of shifted regret decays to zero at a prescribed exponential rate. This optimal decision arises as the solution to a multi-objective optimization problem, which reflects asymptotic behavior properties of the data-generating process. Central to our framework is a rate function that characterizes this behavior via a Laplace principle, thereby generalizing standard concepts from large deviation theory. Our general formulation enables our approach to account for data from uncertain distributions and recovers classical results in data-driven decision making under uncertainty as special cases, including distributionally robust optimization. Moreover, our method enables decision-makers to systematically balance a desired rate of asymptotic risk decay against a potential loss in statistical consistency of the resulting data-driven decision. We demonstrate the effectiveness of the proposed approach through illustrative examples from operations research, such as the newsvendor problem, under aleatoric uncertainty induced by heterogeneous data sources.

\vspace{0.5cm}
\section{Introduction}\label{sec_introduction}
Decision-making problems under uncertainty can be abstractly formulated as $\min_{x \in \dspa} c(x, \theta)$, where $\dspa$ represents the feasible set of decisions, $\theta \in \Theta$ denotes an unknown (potentially infinite-dimensional) distribution parameter, and $c(x, \theta)$ stands for the risk or cost associated with the decision $x$ under the stochastic model characterized by $\theta$. 

Since the parameter $\theta$ is unknown, the problem cannot be directly solved. Nevertheless, we assume access to training sample paths of a $\sspa$-valued stochastic process $\seq{\proc}{n}$, which is distributed according to a probability measure $\mathbb{P}_{\theta}$. For example, $\proc_n$ might be an empirical expectation or an empirical probability measure comprised of the first $n$ elements of some underlying stochastic process. Data-driven decision-making then aims to find a mapping $\dec\pcolon\sspa\map\dspa$ such that $\dec (\proc_n)$ is a ``good'' decision when evaluated at the true unknown model $\theta$, where ``good'' can be expressed via possibly low \defformat{out-of-sample (or generalization) cost} $c(\dec (\proc_n), \theta)$  or via small enough \textit{regret} $c(\dec (\proc_n), \theta) - \min_{x \in \dspa} c(x, \theta)$. 
A fundamental quantity to assess the asymptotic performance of a decision function $\dec$ is the asymptotic counterpart of regret associated with the process $\seq{\proc}{n}$ called \defformat{consistency gap of $\dec$}, which we define as
\begin{equation}\label{eq_consistency_gap_intro}
   \lim_{n\to\infty} \regret(\dec (\proc_n),\theta) = \lim_{n\to\infty}  c(\dec (\proc_n), \theta) - \min_{x \in \dspa} c(x, \theta),
\end{equation}
assuming that the limit exists.
If \eqref{eq_consistency_gap_intro} equals to zero 
$\mathbb{P}_{\theta}\text{-a.s.}$, the underlying data-driven decision is (strongly) consistent.
Clearly, there are various approaches to constructing a data-driven decision that is consistent. For example, if $c$ is expressed as an expected loss, the classical plug-in estimator (also known as \acf{SAA}), where the unknown expectation is replaced with its empirical counterpart, leads under relatively mild conditions to a consistent data-driven decision $\decplug$ defined via
$$
\decplug(\proc_n)\in\Argmin_{x\in\dspa} c(x,\proc_n),
$$
as established by \cite{ref:Wets-91} and \cite{ref:Shaprio-02}. 
Consistency can also be established for other widely used approaches to data-driven decision making under uncertainty. For instance, in the case of \acf{DRO} problems. 
\ac{DRO} aims to make decisions that are optimal against the worst-case probability distribution within a set of plausible distributions, known as ambiguity set, defined by available information (\cite{ref:Kuhn:DRObook-24}). In \ac{DRO}, consistency holds provided that the radius of the ambiguity set decreases at an appropriate rate~(\cite{Wass17Monh}). In the limit, as the number of data points tends to infinity, the resulting decision then converges to that of the plug-in (or \ac{SAA}) problem.

It is well-known (\cite{smith2006optimizer}) that while the plug-in method is consistent, it exhibits poor out-of-sample performance in the finite-sample regime, a phenomenon referred to as the optimizer's curse. Consequently, relying solely on the consistency gap to assess the quality of a data-driven decision appears less reasonable, suggesting the need for alternative evaluation criteria. 
A natural criterion would be to require the regret of the evaluated data-driven decision $\dec (\proc_n)$ under the true model $\theta$ to decay at a certain, say exponential, rate. However, this is challenging for two reasons: first, since $\theta$ is unknown, this decay would need to be guaranteed for any $\theta \in \Theta$; second, achieving an exponential decay rate of the regret might be inherently impossible for some models. Therefore, we propose 
to search for decision functions $X$ that are characterized as optimal solutions to the following criterion
\begin{equation}
\begin{array}{rl}
    \displaystyle\inf_{\dec,\off} &\{\off( z)\}_{ z\in\sspa}\\  
    \displaystyle\subjectto &\limsup_{n\to\infty} \frac{1}{n} \log \mathbb{E}_\theta \left[\exp( n \, \ell(\regret(\dec(\proc_n),\theta)-\off(\proc_n)))\right]\leq -\rate,\quad \forall \theta\in\Theta,
\end{array}
\label{eq_meta_problem_intro}
\end{equation}
where the optimization variables $\dec\pcolon\sspa\map\dspa$ and $\off\pcolon\sspa\map\real$ are functions with certain regularity conditions that will be specified, $\ell\pcolon \realinf\to\realinf$ is a penalty function and $r>0$ is a fixed parameter.
The constraint in \eqref{eq_meta_problem_intro} corresponds to an inequality for the \defformat{asymptotic entropic risk measure} 
$\phi^{\mathsf{ENT}}_{\theta}$, i.e., to $\phi^{\mathsf{ENT}}_{\theta}[(f(\proc_n))_{n\in\nat}]\leq -r<0$ with $f(z) \defeq \ell(\regret(\dec(z),\theta)-\off(z))$ for $z\in\sspa$. That is, problem~\eqref{eq_meta_problem_intro} admits data-driven decisions $\dec$ and upper confidence bounds $\off$ such that the corresponding asymptotic entropic risk of the penalty of the shifted regret is no larger than $-\rate$. This behavior is enforced for all $\theta \in \Theta$, including the true (unknown) model $\theta$. Our criterion~\eqref{eq_meta_problem_intro} shall identify a tuple of an optimal decision and the smallest upper confidence bound that can achieve the required bound on the asymptotic entropic risk, which controls the decay rate of the corresponding regret. That is also why we call $\rate>0$ the \defformat{convergence rate parameter}.
Although \eqref{eq_meta_problem_intro} is just an instance of a problem class introduced below, it is helpful in developing intuition for our framework. 
The constraint in \eqref{eq_meta_problem_intro} is quite general with respect to the choice of $\ell$. By selecting a degenerate type of a penalty function it reduces to
\begin{equation}\label{eq:constraint:intro:prob:special:case}  
\limsup_{n\to\infty} \frac{1}{n} \log \probsub{\theta}{\regret(\dec(\proc_n),\theta) > \off(\proc_n)} \leq -\rate, \quad \forall \theta\in\Theta.
\end{equation}
In this formulation, the interpretation of $\off$ as a high-probability upper confidence bound becomes particularly clear.

The uncertainty considered so far originates from a lack of knowledge about the true data-generating probability vector~$\theta$. This type of uncertainty can typically be reduced by collecting more data. For example, one may consider a setting in which the underlying stochastic process~$(\proc_n)_{n\in\mathbb{N}}$ converges to $\theta$ as $n \to \infty$. Such uncertainty is commonly referred to as \emph{epistemic uncertainty}. Most of the existing literature on data-driven decision making, including stochastic programming~(\cite{ref:Shapiro-09}) and \ac{DRO} (\cite{ref:Kuhn:DRObook-24}), focuses on this notion of uncertainty.
In contrast, \emph{aleatoric uncertainty} refers to randomness that is intrinsic to the data itself. Unlike epistemic uncertainty, aleatoric uncertainty cannot be eliminated, even with access to infinite data~(\cite{ref:wimmer-23}). This type of uncertainty naturally arises in practical scenarios, such as when data is collected from heterogeneous sources. Suppose each source~$k$ corresponds to a slightly different model~$\theta_k \in B_\theta(R)$, where $B_\theta(R)\subseteq\Theta$ is a ball with radius $R$ centered around a nominal (but unknown) model~$\theta$. We assume that the test model, i.e., the one under which the learned policy is deployed, also lies in $B_\theta(R)$.
In such settings, the constraint in \eqref{eq_meta_problem_intro}, which is formulated as an asymptotic entropic risk under the model~$\theta$, naturally generalizes to the asymptotic version of the so-called \defformat{asymptotic robust entropic risk measure}~\cite[Section 5]{iP21}
\begin{equation}\label{eq:intro:robust:entropic:risk}
    \lim_{n \to \infty} \frac{1}{n} \log \sup_{\utheta \in B_\theta(R)} \expecsub{\utheta}{\exp\left(n \, \ell\left(\regret(\dec(\proc_n),\theta) - \off(\proc_n)\right)\right)} \leq -r.
\end{equation}

In this work, we show that under certain assumptions to be specified, the solution to~\eqref{eq_meta_problem_intro}, where the constraint can be expressed as a general asymptotic risk measure which satisfies an abstract Laplace principle (see~\cite{iP10}), encompassing both the asymptotic and asymptotic robust entropic risk constraints, is remarkably simple and takes the form
\begin{equation}\label{eq:solution:meta:problem:intro}
\begin{aligned}
    \optimal{\dec}(z) &\defeq \argmin_{x \in \dspa} \max_{\theta \in \Theta} \left\{ \regret(x,\theta) - \inv{\ell} \left(\ratefct_\theta(z) - \rate \right) \right\}, \\
    \optimal{\off}(z) &\defeq \min_{x \in \dspa} \max_{\theta \in \Theta} \left\{ \regret(x,\theta) - \inv{\ell} \left(\ratefct_\theta(z) - \rate \right) \right\}, \quad z \in \sspa.
\end{aligned}
\end{equation}
The function $\ratefct\pcolon\Theta\times\sspa\to[0,\infty]$ stands in relationship to the constraint in~\eqref{eq_meta_problem_intro} and reflects the asymptotic behavior of $\seq{\proc}{n}$. In the simple setting where the constraint is expressed by the asymptotic entropic risk measure, as in~\eqref{eq_meta_problem_intro}, the function $\ratefct_\theta(\cdot)$ corresponds to the large deviation rate function associated with the family of distributions of the sequence~$\seq{\proc}{n}$. This implies that the structure of the optimal data-driven decision function $\optimal{\dec}$ and upper confidence bound $\optimal{\off}$ is then induced by the concentration behavior of the underlying data generating process via the corresponding large deviation principle, which is represented by the rate function.
In the robust setting, where the constraint in \eqref{eq_meta_problem_intro} is replaced with \eqref{eq:intro:robust:entropic:risk}, the standard large deviation rate function $\ratefct_\theta(\cdot)$ is replaced with by its robust counterpart $\inf_{\utheta\in B_\theta(R)} \ratefct_{\utheta}(\cdot)$.

\paragraph{Related work.}
The optimality criterion~\eqref{eq_meta_problem_intro} is conceptually closely related to universal composite hypothesis testing problems and the underlying asymptotic optimality concept; see \cite[Section 7.1]{iB1} for a comprehensive treatment. This connection is best exemplified by considering the special case where the constraint in \eqref{eq_meta_problem_intro} reduces to \eqref{eq:constraint:intro:prob:special:case}.
Analogous to hypothesis testing problems, we interpret the unfavorable event $\regret(\dec(\proc_n), \theta) > \off(\proc_n)$, where the regret of our decision exceeds the upper confidence bound, as a type I error. Clearly, one could minimize the corresponding error probability by choosing the upper confidence bound $\off$ to be very large. Our optimality criterion~\eqref{eq_meta_problem_intro} then quantifies a pair of data-driven decision and corresponding upper confidence bound $(\optimal{\dec}, \optimal{\off})$ as optimal if it leads to the desired decay rate $r$ of the type I error probability, while at the same time being the smallest upper confidence bound. Specifically, the pair $(\optimal{\dec}, \optimal{\off})$ must satisfy \eqref{eq:constraint:intro:prob:special:case}, and, in addition, any other pair $(\dec, \off)$ satisfying this property must be such that
\begin{equation*}
\liminf_{n \to \infty} \frac{\off(\proc_n)}{\optimal{\off}(\proc_n)} \geq 1.
\end{equation*}
While the optimal composite hypothesis test \cite[Theorem 7.1.11]{iB1} is defined via the level-set of an underlying large deviation rate function, a highly related structure is present in our setting. In fact, the solution \eqref{eq:solution:meta:problem:intro} for the special penalty leading to problem \eqref{eq_meta_problem_intro} with the constraint \eqref{eq:constraint:intro:prob:special:case} reduces to
\begin{equation}\label{eq:optimal:decision:special:case:regret}
    \optimal{\dec}(z) = \argmin_{x \in \dspa} \sup_{\theta \in \Theta,\ \ratefct_\theta(z)< \rate} \regret(x,\theta),
\end{equation}
where we also optimize over the level-sets of the large deviation rate function, but with respect to the $\theta$-variable. The optimization problem~\eqref{eq:optimal:decision:special:case:regret} can be directly interpreted as a \ac{DRO} problem, highlighting the inherent connection between the solutions to \eqref{eq_meta_problem_intro} and \ac{DRO}.

In this paper, we show that the optimal decision function defined in \eqref{eq:solution:meta:problem:intro}, for appropriately chosen asymptotic risk measures, recovers both the classical \ac{SAA} decision rule (\cite{ref:Shapiro-09}) and the well-known \acf{DRO} rule~(\cite{ref:Delage-10, ref:Kuhn:DRObook-24}) as special cases. More importantly, our framework provides a principled way to interpolate between these two extremes via parametrization of the penalty function, yielding a rich family of data-driven decision rules that trade off between the statistical consistency of \ac{SAA} and the robustness of \ac{DRO}, particularly in low-data regimes. This flexibility enables decision makers to tailor the decision rule to the specific characteristics and requirements of their application.

A similar perspective for data-driven decision making has been considered in \cite{ref:vanParys:fromdata-17, iP1}, where the authors aim to choose a data-driven decision $\dec$ with the goal to achieve an exponential decay rate of the so-called out-of-sample disappointment probability. Written in our notation, the authors consider a finite-dimensional problem of the form \eqref{eq_meta_problem_intro}, where the constraint is expressed as
\begin{align}\label{eq:const:DK}
\limsup_{n \to \infty} \frac{1}{n} \log \probsub{\theta}{c(\dec(\proc_n), \theta) > \off(\proc_n)} \leq -\rate, \quad \forall \theta \in \Theta,
\end{align}
and $\off(\proc_n)$ is referred to as the in-sample cost. The optimal decision to this approach is given by \eqref{eq:optimal:decision:special:case:regret} where $\regret(x,\theta)$ is replaced by $c(x,\theta)$. 
Moreover, existing works primarily focus on epistemic uncertainty. A notable exception is~\cite{ref:parys-24}, which considers data-driven decision making under noisy data, i.e., data affected by aleatoric uncertainty. However, this work also relies on independence assumptions for the underlying data samples.

In the context of \ac{DRO}, the trade-off between robustness (e.g., in the sense of \eqref{eq:const:DK}) and the conservatism of the underlying decision (e.g., quantified via the upper-confidence bound $u$) is typically expressed by the size of the ambiguity set. The size of these ambiguity sets can also be interpreted as the strength of a certain regularization, which critically determines the resulting out-of-sample behavior~(\cite{ref:Soroosh:JMLR-19}).
In most works~(see the comprehensive textbook \cite{ref:Kuhn:DRObook-24}), the size of the ambiguity set is chosen via cross-validation. While this approach often works well in practice, it lacks rigorous statistical guarantees.
There are a few approaches that theoretically study the statistically optimal choice of \ac{DRO} ambiguity set sizes.
\cite{ref:Lam-16} proposes ambiguity sets that offer optimal statistical guarantees in view of the central limit theorem.
Moreover, the asymptotic size of \ac{DRO} ambiguity sets with respect to the sample size is known both for divergence ambiguity sets (\cite{ref:Lam-16,ref:Duchi-21}) as well as for Wasserstein ambiguity sets (\cite{ref:Blanchet_2019,gao2020finitesample}).
In a recent work by~\cite{ref:Blanchet:tutorial-21}, the choice of optimal radius of Wasserstein ambiguity sets, in a \ac{DRO} context, has been analyzed via a specific projection that is defined as a hypothesis test.

Our approach is not only closely related to \ac{DRO}, but also conceptually connected to the principles of robust statistics.
Specifically, it can be shown (\cite{ref:Amine-24}) that the ambiguity sets of \ac{DRO} problems based on relative entropy and total variation distance correspond to confidence sets with the desired properties of the optimal confidence set proposed by \cite{ref:Huber-64}. A related result is proposed by \cite{ref:Sut-Gang-23}, where it is demonstrated that among all confidence intervals that asymptotically guarantee a certain coverage probability (referred to as accuracy), a \ac{DRO}-based confidence interval is minimal in size. 

\paragraph{Contribution.}
We highlight the following main contributions of this paper.
\begin{enumerate}[label=\enumintext]
    \item We propose a framework for optimal data-driven decision making in which the optimality criterion generalizes~\eqref{eq_meta_problem_intro} by replacing the asymptotic entropic risk measure constraint with an arbitrary asymptotic risk measure which satisfies a Laplace principle. This generalized formulation enables the representation of a broad class of asymptotic risk measures under minimal assumptions on the underlying stochastic process. In particular, it includes asymptotic robust risk measures that effectively capture aleatoric uncertainty.

    \item We show that the proposed generalizations of problem~\eqref{eq_meta_problem_intro} admit remarkably simple closed-form solutions whose structure is induced by the statistical properties of the corresponding stochastic process $\seq{\proc}{n}$. Specifically, assuming that the distribution of $\seq{\proc}{n}$ is such that a Laplace principle for a given asymptotic risk measure can be established, the solution is constructed from the underlying rate function, thereby incorporating information about asymptotic behavior of $\seq{\proc}{n}$. Moreover, we demonstrate universality of this optimal solution structure, as it can be shown to solve other modifications of the generalized problem~\eqref{eq_meta_problem_intro} as well.
    
    \item We analyze how the choice of the penalty function~$\ell$ in the generalized formulations of~\eqref{eq_meta_problem_intro} influences the consistency, or the resulting consistency gap, of the optimal data-driven decision. These insights offer practical design guidelines for selecting the penalty function $\ell$.
    In particular, these parameters allow for a continuous interpolation between the classical and well-studied plug-in and \ac{DRO} decision rules.

    \item We derive explicit optimal data-driven decisions for various settings, including finite-state and continuous i.i.d.\ processes with unknown means, as well as finite-state Markov chains. As a case study, we apply our framework to classical OR problems such as the newsvendor problem, both in the standard setting with unknown demand and a more complex variant with additional aleatoric uncertainty due to random demand distributions. These examples, along with supporting numerical results, illustrate the expressiveness and broad applicability of our framework.

\end{enumerate}

\paragraph{Structure.}
Section~\ref{sec_problem_statement} formalizes the problem setting and introduces the generalized form of the optimality criterion~\eqref{eq_meta_problem_intro} outlined in the introduction. It also presents the corresponding solutions and discusses the associated consistency gap.  
Section~\ref{sec_other_opt_criteria} presents the most general formulation of our data-driven optimization framework and its modifications, demonstrating the robustness of the solution structure to changes in the problem setup.
Section~\ref{sec_examples} provides concrete modeling examples and presents numerical studies of two classical Operations Research problems: the newsvendor problem and optimal portfolio selection. 
Finally, Section~\ref{sec_proofs} contains the proofs of the main results from the previous sections.  
Additional technical details and auxiliary results are provided in the \nameref{sec_appendix}.

\paragraph{Notation.}
Let $f\pcolon Y\map\realinf$ be a function taking values in the extended real line.
The \defformat{(effective) domain} of $f$, denoted $ \dom_{Y} f$, is then defined as $\dom_{Y} f\defeq\set{y\in Y}{f(y)<\infty}$. 
The set of minimizers of $f$ is denoted by $\Argmin_{y\in Y} f(y)$ and if this set is a singleton, its unique element is denoted by $\argmin_{y\in Y} f(y)$. Next, $\indfctset{A}$ is the indicator function on some set $A\subseteq Y$.
Let $Y_1,Y_2$ be some topological spaces. We write $\ctsset(Y_1,Y_2)$ for the set of continuous functions $f\pcolon Y_1\to Y_2$, which is abbreviated to $\ctsset(Y_1)$ if $Y_2=\real$. Then, $\ctsset_b$ and 
$\ctsset_{ub}$ denote the subsets of bounded and from above bounded continuous functions respectively.
The interior of some set $A\subseteq Y_1$ is denoted by $\interior{A}$ and $\overline{A}$ stands for its closure.
Moreover, $\mathcal{M}_1(\salg)$ is the space of probability measures on some measurable space $(\Omega,\salg)$ and if $\Omega$ is a topological space, $\mathcal{B}(\Omega)$ stands for the Borel $\sigma$-algebra on $\Omega$.
Finally, $e_1,\ldots,e_d$ are the canonical basis vectors of $\real^d$.

\section{Problem Statement and Main Results}\label{sec_problem_statement}

Let $\seq{\proc}{n}$ be a sequence of $\sspa$-valued random variables defined on a probability space $(\pspa,\salg,\bbP_\theta)$, where $\sspa$ is some non-empty metric state space and $\theta$ is a distribution parameter in some non-empty metric space $\Theta$. Further, let $\dspa$ be another non-empty metric space called the \defformat{decision space}. We assume both $\Theta$ and $\dspa$ to be compact and convex.
Let $c\pcolon \dspa\times \Theta\to\real$ be a continuous \defformat{cost function} so that $x\mapsto c(x,\theta)$  is strictly convex for all $\theta\in\Theta$.
The key problem of interest is then
\begin{equation} \label{eq_overall_optp}
\min_{x \in \dspa} c(x, \theta),
\end{equation}
where $\theta\in\Theta$ is an unknown model parameter. In other words, without the knowledge of $\theta$ and based on the inference of the underlying process $\seq{\proc}{n}$ we look for a cost function minimizing decision $x\in\dspa$, that best approximates \eqref{eq_overall_optp}.
A popular choice of cost $c$ would be described by an expected loss, i.e., $c(x,\theta) = \mathbb{E}_\theta[h(x,\xi)]$, where $\xi$ is an $\Xi$-valued random variable and $h:\mathcal{X}\times\Xi\to\real$ some loss function.

We define a \defformat{data-driven decision (or decision function)} as a continuous function $\dec\pcolon \sspa\to\dspa$, which is to be evaluated at some step of the stochastic process $\seq{\proc}{n}$.
In the special case where $\proc_n \to \theta$ as \(n \to \infty\), i.e., $\proc_n$ is a consistent distribution estimator, a widely used approach in the literature is the \defformat{plug-in estimator} resulting in the \defformat{plug-in decision} defined via
\begin{equation}\label{eq_plug}
    \decplug(\proc_n) \defeq \argmin_{x\in\dspa} c(x,\proc_n).
\end{equation}
When the cost $c$ represents an expected loss, this formulation aligns with the well-known \acf{SAA} method.
Our goal is to find $\dec$ systematically so that $\dec(\proc_n)$ approximates the minimizer in \eqref{eq_overall_optp} ``effectively'', which will be formalized in the following.

To formulate a problem of the type~\eqref{eq_meta_problem_intro}, it is essential to understand how the transformed stochastic process $(\ell(\regret(\dec(\proc_n),\theta) - \off(\proc_n)))_{n \in \mathbb{N}}$, which is influenced by the decision function $\dec$, concentrates. This process captures the penalized shifted regret associated with the evaluated data-driven decision $\dec(\proc_n)$.
We employ asymptotic risk measures, such as the asymptotic entropic risk measure (see~\eqref{eq_meta_problem_intro}), to characterize and control the concentration behavior of the transformed process. 
Problem~\eqref{eq_meta_problem_intro} defines the class of admissible decision functions $\dec$ as those for which the asymptotic risk of the associated transformed stochastic process $(\ell(\regret(\dec(\proc_n),\theta) - \off(\proc_n)))_{n \in \mathbb{N}}$ decays exponentially at a prescribed rate $r > 0$.
Motivated by the theory of large deviations, we introduce two mathematical objects that together encapsulate such information about the (transformed) stochastic process $\seq{\proc}{n}$.

First, let $\phi_\theta$ be a functional called the \defformat{asymptotic risk measure of $\seq{\proc}{n}$ under $\theta$}, which maps any continuous transformation $(f(\proc_n))_{n\in\nat}$ of the stochastic process $\seq{\proc}{n}$ to the extended real line $\realinf$, i.e., 
\begin{align*}
    (f(\proc_n))_{n\in\nat} \mapsto \phisub{\theta}{(f(\proc_n))_{n\in\nat}}\quad\text{for }f\in\ctsset(\sspa).
\end{align*}
A widely used instance of such an asymptotic risk measure is the mentioned asymptotic entropic risk measure represented by 
\begin{equation*}
    (f(\proc_n))_{n\in\nat}\mapsto\limsup_{n\to\infty}\frac{1}{n}\log\bbE_\theta\left[\exp\left(n f(\proc_n)\right)\right]
\end{equation*}
for $f\in\ctsset_{b}(\sspa)$ and we shall see more examples as we go along.
In the generalized setting, the functional $\phi_\theta$ then replaces the asymptotic entropic risk measure in problem \eqref{eq_meta_problem_intro} and thereby defines the feasibility characterizing probabilistic constraint.
Such an abstract construction might appear more intuitive, when compared to the ambivalence between the study of the expected values of arbitrary bounded continuous transformations of some random variable and the study of its distribution. Hence, this is in parallel to studying evaluations of $\phi_\theta$ on arbitrary continuous transformations of a stochastic process in order to encapsulate some of its asymptotic behavior properties.
In other words, in our analogy $\phi_\theta$ replaces the expectation functional $\bbE_\theta$ as the focus shifts from distribution of a random variable to concentration properties of a stochastic process.

Second, the asymptotic behavior of the stochastic process $\seq{\proc}{n}$ described by the asymptotic risk measure $\phi_\theta$ shall be exactly such that there exists a so called \defformat{rate function} $$\ratefct\pcolon \Theta\times\sspa\map[0,\infty],\quad (\theta,z)\mapsto \ratefct_\theta(z),$$ which together with $\phi_\theta$ constitutes a \acf{LP}
\begin{equation}\label{eq_LP}
    \phisub{\theta}{(f(\proc_n))_{n\in\nat}}=\supb{z\in\mathcal{Z}}{f(z)-\mathcal{I}_\theta(z)},\quad \forall f\in\ctsset(\sspa),\ \forall\theta\in\Theta.
\end{equation}
Then, we say that the pair $(\phi_\theta,\ratefct_\theta)$ satisfies the \ac{LP}~\eqref{eq_LP}. Consequently, all the objects $\phi_\theta$, $\ratefct_\theta$ and $\seq{\proc}{n}$ stand in a mutual relationship, which is encoded by \eqref{eq_LP}.
Moreover, with \eqref{eq_LP} at hand, one can interpret the shape of $\ratefct_\theta$ in a direct relationship to the concentration properties of $\seq{\proc}{n}$
and with the optimal decision formula \eqref{eq:solution:meta:problem:intro} in mind also to the consistency properties of such decision.
This also provides an insight into the inference of information about $\theta$ via observations of $\seq{\proc}{n}$.

Before proceeding with the general setup, we consider the specific example where $\phi_\theta$ is chosen as the asymptotic entropic risk measure, in order to develop a better understanding of~\eqref{eq_LP}.

\begin{example}[LDP - asymptotic entropic risk measure]\label{ex_standardLDP_ARE}
    In order to arrive at the \ac{LP}~\eqref{eq_LP} with $\phi_\theta$ being the asymptotic entropic risk measure, 
    we take a look into the theory of large deviations--a literature area, which largely inspires our nomenclature. Let $\theta\in\Theta$ and denote the $\theta$-dependent family of distributions (or laws) associated with $\seq{\proc}{n}$ by $\seq{\law^{\theta}}{n}$, i.e., $\law^{\theta}_n(A) \defeq \probsub{\theta}{\proc_n\in A}$ for any $A\in\mathcal{B}(\sspa)$ and $n\in\nat$. Assume that $\seq{\law^{\theta}}{n}$ satisfies a so called \acf{LDP} with the rate function $\ratefct_\theta$, i.e., 
    \begin{equation}\label{eq_standardLDP}
        -\inf_{z\in \interior{A}} \ratefct_\theta(z) \leq \liminf_{n\converge\infty}\frac{1}{n}\log \law_n^{\theta}(A)\leq \limsup_{n\converge\infty}\frac{1}{n}\log \law_n^{\theta}(A)\leq-\inf_{z\in \overline{A}} \ratefct_\theta(z)
    \end{equation}
    hold for all Borel sets $A\subseteq\sspa$. If $\ratefct$ is such that $\ratefct_\theta\pcolon\sspa\map[0,\infty]$ is a \textit{good rate function}\footnote{
    If $\ratefct_\theta$ is an \ac{lsc} function with compact \defformat{level sets} $\set{z\in\sspa}{\ratefct_\theta(z)\leq\alpha}$, $\alpha\in[0,\infty)$, then it is called a \defformat{good rate function} according to the standard large deviation literature \cite[Section 1.2]{iB1}.
    },
    then Varadhan's Lemma \cite[Theorem 4.3.1]{iB1} guarantees the \ac{LP}
    \begin{equation}\label{eq_standardVaradhan}
        \phiupsub{ENT}{\theta}{(f(\proc_n))_{n\in\nat}}\defeq
        \lim_{n\to\infty}\frac{1}{n}\log\bbE_\theta\left[\exp\left(n f(\proc_n)\right)\right]=\supb{z\in\mathcal{Z}}{f(z)-\mathcal{I}_\theta(z)}
    \end{equation}
    for all $f\in\ctsset_{ub}(\sspa)$ with the asymptotic entropic risk measure\footnote{In standard large deviation theory the so called \defformat{asymptotic relative entropy}.} $\phi^{\mathsf{ENT}}_{\theta}$. Notice, for $\phi^{\mathsf{ENT}}_{\theta}$ to truly be an asymptotic risk measure of $\seq{\proc}{n}$ under $\theta$ with the \ac{LP}~\eqref{eq_LP}, \eqref{eq_standardVaradhan} needs to hold for any $f\in\ctsset(\sspa)$.
    We ensure this by extending $\phi^{\mathsf{ENT}}_{\theta}$ to all continuous transformations of the underlying process $\seq{\proc}{n}$, i.e., $(f(\proc_n))_{n\in\nat}$ with $f\in\ctsset(\sspa)$, via
    \begin{equation}\label{eq_extendedARE}
        \phiupsub{ENT}{\theta}{(f(\proc_n))_{n\in\nat}}\defeq\sup_{m\in\real}\limsup_{n\to\infty}\frac{1}{n}\log\expecsub{\theta}{\exp(n \min\{f(\proc_n),m\})}.
    \end{equation}
    It can be shown that, with such an extension, the \ac{LP}~\eqref{eq_LP} indeed holds for any $f\in\ctsset(\sspa)$  \cite[Theorem 3.5, Proposition 7.3]{iP11}.
\end{example}
Our terminology differs slightly from the standard definition of an \ac{LDP} and a rate function \cite[p.4]{iB1}. For instance, the parametrization of the rate function $\ratefct$ in $\theta$ is usually not considered, so that for any fixed $\theta$ the function $\ratefct_\theta\pcolon\sspa\to[0,\infty]$ is the actual object commonly known as rate function.
An example of an \ac{LDP} follows.

\begin{example}[Finite state i.i.d. process with expected loss]\label{ex_main_base}
    Let $\Delta_d\defeq\{x\in\real_{\geq0}^d\mid\sum_{i=1}^d x_i=1\}$ be the \defformat{probability simplex} in $\real^d$ for some $d\in\nat$ and consider $\Theta\subseteq\Delta_d$. 
    Let $\seq{\xi}{k}$ be an i.i.d. finite state stochastic process on the probability space $(\pspa,\salg,\bbP_\theta)$ such that $\xi_k\pcolon \pspa\map
    \{1,\ldots,d\}$ and $ \probsub{\theta}{\xi_k=i}=\theta_i$ for $i=1,\ldots,d$.
    Let $\proc_n$ be the empirical measure associated with $(\xi_1,\ldots,\xi_n)$, i.e., $\proc_n = \sum_{i=1}^d \left(\frac{1}{n}\sum_{k=1}^{n}\indfctset{\xi_k=i}\right)e_i$ and $\sspa=\Delta_d$. Further, we define the \defformat{relative entropy} $\ent\pcolon\Theta\times\sspa\subseteq\Delta_d\times\Delta_d\to[0,\infty]$ by
    \begin{equation}\label{eq_def_relative_entropy}
        \ent_\theta ( z)\defeq\begin{cases}
    \sum_{i=1}^d z_i \log(z_i/\theta_i)& \text{ if } \theta_i =0\Rightarrow z_i=0\\
    \infty&\text{ otherwise},
    \end{cases}
    \end{equation}
    where the standard conventions $0\log0=0$ and $0\log(0/0)=0$ apply.
    Then, by Sanov's theorem \cite[Theorem 2.1.10]{iB1} the sequence $\seq{\law^{\theta}}{n}$, as defined in \Cref{ex_standardLDP_ARE}, satisfies a standard \ac{LDP} with the good rate function $\ent_\theta$. 
    A suitable choice of cost $c$ would be an expected loss, i.e., 
    $c(x,\theta) =\sum_{i=1}^d h(x,i) \theta_i$,
    where $h\pcolon \dspa\times\{1,\dots,d\}\to\real$ is continuous and strictly convex in $x$.
\end{example}

\begin{example}[Law of large numbers] \label{ex_LP_via_LLN}
    A natural way to demonstrate the ability of the asymptotic risk measure to capture the concentration behavior of the process $\seq{\proc}{n}$, and that an \ac{LDP} from standard large deviation theory is merely one of many options of establishing a Laplace principle, is to observe the  case when $\seq{\proc}{n}$ satisfies a \ac{LLN}, i.e., $\proc_n\to\theta$ as $n\to\infty$ holds $\bbP_\theta$-a.s. For instance the empirical measure process in the previous \Cref{ex_main_base} has this property. Then, the simplest way of constructing an \ac{LP}~\eqref{eq_LP} is to consider the asymptotic risk measure
    \begin{equation*}
        \phiupsub{LLN}{\theta}{(f(\proc_n))_{n\in\nat}}\defeq f(\theta) = \lim_{n\to\infty} f(\proc_n),\quad f\in\ctsset(\sspa), \ \bbP_\theta\text{-a.s.},
    \end{equation*} 
    because, provided $\Theta\subseteq\sspa$, it is obvious that with the rate function 
    \begin{equation*}
     \ratefct^{\mathsf{LLN}}_\theta(z)   \defeq\begin{cases}
   0& \text{ if } z = \theta\\
    \infty&\text{ otherwise},
    \end{cases} \qquad (\theta,z)\in\Theta\times\sspa,
    \end{equation*}
    we get
    \begin{equation*}
        \phiupsub{LLN}{\theta}{(f(\proc_n))_{n\in\nat}} = f(\theta) = \sup_{z\in\sspa}\{f(z)-\ratefct^{\mathsf{LLN}}_\theta(z)\}, \quad\forall f\in\ctsset(\sspa), \ \forall \theta\in\Theta.\qedhere
    \end{equation*}
\end{example}

We now return to the central optimization problem of the form \eqref{eq_meta_problem_intro}, which constitutes the core of this work.
Let $\ell\pcolon \real\map\real$ be a \defformat{penalty function}, which is an increasing continuous and bijective function with $\ell(0)=0$. Hence, it holds $\lim_{y\to\pm\infty} \ell(y)=\pm\infty$ and we set $\ell(\infty) \defeq \infty$ and $\ell(-\infty) \defeq -\infty$. More important, however, is the inverse of the penalty function, since the optimal decision function, as already displayed in \eqref{eq:solution:meta:problem:intro}, will be defined in terms of $\inv{\ell}$ only. 
For example, we define the penalty function $\ell_{\alpha,\beta}$ via its inverse 
\begin{equation}\label{eq_def_penalty_alphabeta}
    \inv{\ell}_{\alpha,\beta}(y)\defeq \frac{1}{\beta}(\exp(\beta y) -1) + \alpha y ,\quad y\in\real,
\end{equation}
with some fixed parameters $\alpha,\beta>0$, see \Cref{fig_plot_penalty_alphabeta_beta}. 

\subsection{Problem statement} \label{subsec_problem_statement_}
A natural quantity to assess the quality of a decision $x \in \dspa$ for problem~\eqref{eq_overall_optp} is its corresponding \defformat{regret} under the distribution specified by $\theta \in \Theta$, defined as
\begin{equation*}
    \regret(x, \theta) \defeq c(x, \theta) - \min_{y \in \dspa} c(y, \theta).
\end{equation*}
We formulate the main optimization problem of interest in this work as
\begin{equation}\label{eq_optp_main} 
\begin{array}{cl}  \displaystyle\inf_{\dec\in\decset, \ \off\in\uset}& \{\off( z)\}_{ z\in\sspa}\\  
    \displaystyle \ \subjectto &\phisub{\theta}{\left(\ell\left(\regret(\dec(\proc_n),\theta)-\off(\proc_n)\right)\right)_{n\in\nat}}\leq -\rate,\quad \forall \theta\in\Theta,
\end{array} 
\end{equation}
where $\rate > 0$ is a fixed parameter, which can in certain instances be interpreted as a constant that controls the asymptotic decay rate of the unfavorable event 
$\{\regret(\dec(\proc_n), \theta) > \off(\proc_n)\}$ and is therefore called the \defformat{convergence rate} parameter. To see this, consider the special case described in Example~\ref{ex_standardLDP_ARE}, where $\phi_\theta$ represents the asymptotic entropic risk measure and suppose the penalty function is chosen such that the constraint in \eqref{eq_optp_main} is approximately characterized by  
\begin{equation*}
    \limsup_{n\to\infty} \frac{1}{n} \log \probsub{\theta}{\regret(\dec(\proc_n),\theta) > \off(\proc_n)} \leq -\rate, \quad \forall \theta\in\Theta,\footnote{This can be achieved with $\inv{\ell}_{\alpha,\beta}$ for some large $\beta\gg0$ and small $\alpha\approx0$, see \Cref{rem_discussion_problem_solution}.}
\end{equation*}
which can be intuitively reformulated as $\probsub{\theta}{\regret(\dec(\proc_n),\theta) > \off(\proc_n)} \leq \exp(-nr)$ for $n$ large enough. Moreover, the feasibility constraint in \eqref{eq_optp_main} allows only for such continuous functions $\dec\pcolon\sspa\to\dspa$ and $\off\pcolon\sspa\to\real$, for which the asymptotic risk measure $\phi_\theta$ 
evaluated on the continuous transformation $\ell\left(\regret(\dec(\proc_n),\theta)-\off(\proc_n)\right)$, $n\in\nat$, is bounded by the negative of the convergence rate $\rate$ across all $\theta\in\Theta$. Again in the context of \Cref{ex_standardLDP_ARE}, 
this constraint can be interpreted as an upper bound on the asymptotic entropic risk of the random process $(\ell\left(\regret(\dec(\proc_n),\theta)-\off(\proc_n)\right))_{n\in\nat}$, which reminds us of the introductory optimization problem \eqref{eq_meta_problem_intro}. Hence, problem \eqref{eq_optp_main} identifies a data-driven decision $\dec$ and the smallest corresponding upper confidence bound $\off$ such that the corresponding asymptotic risk of the penalty of the shifted regret is no larger than $-\rate$.

As \eqref{eq_optp_main} is a multi-objective problem, we shall properly define the solution characterizing concepts.
We call the pair consisting of some $\dec\pcolon \sspa\map\dspa$ and $\off\pcolon \sspa\map\real$ \defformat{feasible} in \eqref{eq_optp_main}, if $(\dec, \off)\in\decset\times\uset$ and if it satisfies the constraint therein, i.e., if $$\phisub{\theta}{\left(\ell\left(\regret(\dec(\proc_n),\theta)-\off(\proc_n)\right)\right)_{n\in\nat}}\leq -\rate,\quad\forall\theta\in\Theta.$$ A feasible pair $(\dec, \off)$ is \defformat{optimal} in \eqref{eq_optp_main}, if for any other feasible pair $(\dec^\prime, \off^\prime)$ it holds $$\off(z)\leq\off^\prime(z),\quad \forall z\in\sspa.$$
In other words, the notion of optimality in the context of the multi-objective optimization problem~\eqref{eq_optp_main} refers to \defformat{Pareto dominant solutions}, i.e.,
feasible solutions that Pareto dominate all other feasible solutions. A weaker notion of optimality is the so-called \defformat{Pareto optimal solutions}, i.e., feasible solutions that are not Pareto dominated by any other feasible solution. While typically multi-objective optimization problems only admit Pareto optimal solutions, we show that our problem~\eqref{eq_optp_main} admits a Pareto dominant solution that can be expressed in closed form.

\subsection{Optimal decision} 
In the following, we show that under suitable assumptions an optimal solution to problem~\eqref{eq_optp_main} is given by
\begin{align} \label{eq_solution_star}
\begin{split}
    \optimal{\dec}(z) &\defeq \argmin_{x\in\dspa} \maxb{\theta \in \Theta}{ \regret(x,\theta) - \inv{\ell}\left(\ratefct_\theta(z) - \rate\right)} \quad\text{and}\\
    \optimal{\off}(z) &\defeq \min_{x\in\dspa}\maxb{\theta \in \Theta}{ \regret(x,\theta) - \inv{\ell}\left(\ratefct_\theta(z) - \rate\right)}, \quad z \in \sspa.
\end{split}
\end{align}
We also use the notation
\begin{equation}\label{eq_G}
    \supfct(x,z)\defeq \maxb{\theta \in \Theta}{ \regret(x,\theta) - \inv{\ell}\left(\ratefct_\theta(z) - \rate\right)}, \quad (x,z)\in\dspa\times\sspa,
\end{equation}
so that we can write $\optimal{\dec}(z) = \argmin_{x\in\dspa}\supfct(x,z)$ and $\optimal{\off}(z) = \min_{x\in\dspa}\supfct(x,z)$ for any $z \in \sspa$.
The optimal decision function $\optimal{\dec}$ can be interpreted as a regularized robust decision, as it aims to minimize a regularized regret under a worst-case model $\theta$, where the regularization term $\inv{\ell}\left(\ratefct_\theta(\cdot) - \rate\right)$ depends on the concentration properties of the underlying stochastic process via the rate function $\ratefct$.

To show that the pair $(\optimal{\dec},\optimal{\off})$ is indeed an optimal solution to \eqref{eq_optp_main}, we first need to introduce assumptions on the rate function $\ratefct$, which are mainly there to ensure continuity of $\supfct\pcolon \dspa\times\sspa\map\real$ and $\optimal{\dec}\pcolon \sspa\to\dspa$, as well as to guarantee that the respective optima are truly attained, and that $\optimal{\dec}(z)$ is uniquely defined for all $z\in\sspa$.

\begin{assumption}[Regularity of {$\ratefct\pcolon \Theta\times\sspa\map[0,\infty]$}]\label{ass_ratefct}\ Assume the following.
    \begin{enumerate}[label=\enumassumption]
        \item For every $ z \in \sspa$ there exists $\theta\in\Theta$ such that $(\theta,z)\in\dom_{\Theta\times\sspa} \ratefct$.
        \label{item_ass_ratefct_dom}
        
        \item The rate function $\ratefct$ is \ac{lsc} on $\Theta\times\sspa$. 
        \label{item_ass_ratefct_lsc}
        
        \item For every $(\theta,z)\in\dom_{\Theta\times\sspa} \ratefct$ and all $\seq{z}{n}\subseteq\sspa$ with $z_n\to z$ there exists $\seq{\theta}{n}\subseteq\Theta$ such that $\theta_n\to\theta$ and $ \ratefct_{\theta_n}(z_n)\to\ratefct_\theta(z)$ as $n\to\infty$. 
        \label{item_ass_ratefct_edgects}
    \end{enumerate}
\end{assumption}
By recalling the structure of the proposed data-driven decision function $\optimal{\dec}$ in \eqref{eq_solution_star}, it becomes apparent, that the continuity of $\supfct$ is the key for establishing the continuity of $\optimal{\dec}$ itself.
\Cref{ass_ratefct}~\ref{item_ass_ratefct_dom} contributes to this relationship by ensuring $\supfct>-\infty$, as $\ratefct$ is the only relevant object in \eqref{eq_G} possibly counteracting such required behavior. Such a decisive role of $\ratefct$ applies also to the smoothness properties of $\supfct$. 
If we were to assume the continuity of $\ratefct$ on $\Theta\times\sspa$, the continuity of $\supfct$ would follow immediately from the basic versions of Berge's maximum theorem (\Cref{thm_BergeMaximum}). However, this would exclude relevant examples such as the relative entropy (\Cref{ex_main_base}) and hence, we use slightly weaker conditions \ref{item_ass_ratefct_lsc}~\&~\ref{item_ass_ratefct_edgects} in \Cref{ass_ratefct} instead. 
The lower semi-continuity of $\ratefct_\theta\pcolon\sspa\to[0,\infty]$ for any $\theta\in\Theta$ is one of the fundamental characteristics of a rate function in the standard large deviation theory, where the parameter $\theta$ is not explicitly considered to be a variable input of a rate function at all. Rather naturally, in \ref{item_ass_ratefct_lsc} we extend this typical \ac{lsc} assumption to the joint input $(\theta,z)$. In \ref{item_ass_ratefct_edgects} we provide a condition, which is weaker than requiring $\ratefct$ to be continuous, so that it allows $\ratefct$ to take the value $\infty$ and at the same time helps to resolve smoothness behavior on the boundary of $\dom_{\Theta\times\sspa}\ratefct$. With this the continuity of $\supfct$ can be established.

Building up on the smoothness properties of the optimal decision $\optimal{\dec}$, the following theorem states the feasibility and optimality of the pair $(\optimal{\dec},\optimal{\off})$ in the optimization problem~\eqref{eq_optp_main}.
\begin{thm}[Feasibility and optimality of $(\optimal{\dec},\optimal{\off})$]\label{thm_feas_dom_main} 
    Let the pair $(\phi_\theta,\ratefct_\theta)$ satisfy the \ac{LP}~\eqref{eq_LP} and let \Cref{ass_ratefct} hold. Then, the following is true.
    \begin{enumerate}[label=\enumthm]
        \item \label{item_thm_feas_dom_main_feas} The candidate pair $(\optimal{\dec}, \optimal{\off})$  defined in \eqref{eq_solution_star} is a feasible solution to \eqref{eq_optp_main}, i.e., $(\optimal{\dec}, \optimal{\off})\in\decset\times\uset$ and 
        \begin{equation*}
            \phisub{\theta}{\left(\ell\left(\regret(\optimal{\dec}(\proc_n),\theta)-\optimal{\off}(\proc_n)\right)\right)_{n\in\nat}}\leq -\rate,\quad \forall \theta\in\Theta.
        \end{equation*}
        \item \label{item_thm_feas_dom_main_dom} The pair $(\optimal{\dec},\optimal{\off})$ is \defformat{optimal} in \eqref{eq_optp_main}, i.e., for any other feasible pair $(\dec, \off)$ it holds $$\optimal{\off}(z)\leq\off(z),\quad \forall z\in\sspa.$$
    \end{enumerate}
\end{thm}
\Cref{thm_feas_dom_main} holds as a special case of statements about a more general formulation of the optimization problem \eqref{eq_optp_main}, which are discussed in \Cref{sec_other_opt_criteria} and proven in \Cref{sec_proofs}. There, we will see that the regret function in the problem formulation \eqref{eq_optp_main} and solution \eqref{eq_solution_star} might be replaced by the cost function $c$ 
or any other function sharing the minimal requirements imposed on $c$. Such adjustments can lead to some (numerical) advantages. 

\begin{example}[Plug-in decision as optimal solution in LLN context]
\label{ex_plugin_via_LLN}
    We consider the optimal solution according to \Cref{thm_feas_dom_main} in the context of the \ac{LLN} based Laplace principle from \Cref{ex_LP_via_LLN}, i.e., we assume $\proc_n\to\theta$ as $n\to\infty$ holds $\bbP_\theta$-a.s. and we work with the asymptotic risk measure $\phi^{\mathsf{LLN}}_\theta$ and the degenerate rate function $\ratefct^{\mathsf{LLN}}_\theta$ defined therein. It turns out, that the optimal decision $\optimal{\dec}$ given by \eqref{eq_solution_star} then collapses to the well known plug-in decision \eqref{eq_plug}. Indeed, 
    for any $(x,z)\in\dspa\times\sspa$ we have
    \begin{align*}
        \supfct(x,z)&= \maxb{\theta \in \Theta}{ \regret(x,\theta) - \inv{\ell}\left(\ratefct^{\mathsf{LLN}}_\theta(z) - \rate\right)} \\
        &= \regret(x,z) - \inv{\ell}\left(\ratefct^{\mathsf{LLN}}_z(z) - \rate\right) = \regret(x,z) - \inv{\ell}(-\rate),
    \end{align*}
    so that 
    \begin{align*}
        \optimal{\dec}(z)=\argmin_{x\in\dspa}\supfct(x,z) = \argmin_{x\in\dspa} \regret(x,z) - \inv{\ell}(-\rate) = \argmin_{x\in\dspa} c(x,z) = \decplug(z).
    \end{align*}
    Furthermore, it is easy to verify, that \Cref{ass_ratefct} is fulfilled with $\sspa\subseteq\Theta$. As $\Theta\subseteq\sspa$ was necessary for $(\phi^{\mathsf{LLN}}_\theta,\ratefct^{\mathsf{LLN}}_\theta)$ to satisfy the \ac{LP}~\eqref{eq_LP}, the condition $\Theta=\sspa$ ensures \Cref{thm_feas_dom_main} to verify $\decplug$ as the optimal decision.
    In other words, by using the fact that the stochastic process $\seq{\proc}{n}$ converges almost surely to $\theta$ without any other distributional insight, we proceed with $\phi^{\mathsf{LLN}}$ and $\ratefct^{\mathsf{LLN}}$, ending up with the well-known classical plug-in decision function $\decplug$.

    In the context of the finite Sanov setting in \Cref{ex_standardLDP_ARE} with the empirical measure process $\seq{\proc}{n}$, such \ac{LLN} approach would be admissible.
    Alternatively, as this example goes beyond the \ac{LLN} assumption, the decision-maker may consider a different asymptotic risk measure. Namely the asymptotic entropic risk measure $\phi^{\mathsf{ENT}}_\theta$, which in our setting satisfies the \ac{LP}~\eqref{eq_LP} with the relative entropy $\ent$ as the rate function, as justified in \Cref{ex_standardLDP_ARE}. According to \Cref{thm_feas_dom_main}, the corresponding optimal data-driven decision is given by
\begin{equation*}
    \dec^{\mathsf{ENT}}(z) \defeq \argmin_{x \in \dspa} \max_{\theta \in \Theta} \left\{ \regret(x, \theta) - \inv{\ell}\left(\ent_\theta(z) - \rate\right) \right\},\quad z\in\sspa.
\end{equation*}

When comparing the decision functions $\decplug$ and $\dec^{\mathsf{ENT}}$, one observes that the latter exploits the shape of the rate function $\ent$ to provide more distribution-specific information about the concentration behavior—particularly the rate of decay for rare or extreme events. 
However, this additional information entails a trade-off with the statistical consistency of the plug-in decision, which is reflected in the degenerate form of the rate function $\ratefct^{\mathsf{LLN}}$. In the degenerate \ac{LLN} case, the convergence rate $\rate$ and the specific choice of penalty function $\ell$ become obsolete; likewise, the feasibility condition in~\eqref{eq_optp_main} becomes less meaningful. 
This illustrates how the choice of the asymptotic risk measure is intrinsically linked to the assumed level of distributional uncertainty imposed on the sequence $(\proc_n)_{n \in \mathbb{N}}$.
\end{example}

\begin{rem}[Decision controlling the out-of-sample cost]\label{rem_g_is_c}
An alternative optimal decision function to 
$\optimal{\dec}$, as defined in \eqref{eq_solution_star}, can be derived by considering a problem of the form \eqref{eq_optp_main} with $\regret$ replaced by the cost function $c$. That is, we consider 
\begin{equation}\label{eq_optp_main_c} 
\begin{array}{cl}  \displaystyle\inf_{\dec\in\decset, \ \off\in\uset}& \{\off( z)\}_{ z\in\sspa}\\  
    \displaystyle \ \subjectto &\phisub{\theta}{\left(\ell\left(c(\dec(\proc_n),\theta)-\off(\proc_n)\right)\right)_{n\in\nat}}\leq -\rate,\quad \forall \theta\in\Theta,
\end{array} 
\end{equation}
where the corresponding optimal solutions are then given by  
\begin{align} \label{eq_solution_star_c}
\begin{split}
    \optimal{\dec}(z) &= \argmin_{x\in\dspa} \maxb{\theta \in \Theta}{ c(x,\theta) - \inv{\ell}\left(\ratefct_\theta(z) - \rate\right)} \quad\text{and}\\
    \optimal{\off}(z) &= \min_{x\in\dspa}\maxb{\theta \in \Theta}{ c(x,\theta) - \inv{\ell}\left(\ratefct_\theta(z) - \rate\right)}, \quad z \in \sspa.
\end{split}
\end{align}
\Cref{thm_feas_dom_main} applies accordingly as we show below. While \eqref{eq_optp_main} characterizes the smallest upper-confidence bound $\off$ that controls the regret of the corresponding decision $\dec$, the new formulation \eqref{eq_optp_main_c} does control the \defformat{out-of-sample cost} $c(\dec(\cdot),\theta)$ corresponding to the decision $\dec$. 
This idea, discussed also in the following \Cref{rem_discussion_problem_solution}, is also studied in \cite{ref:vanParys:fromdata-17,iP1}, albeit in a slightly different formulation.
Moreover, it can be seen that the evaluated data-driven decisions $\optimal{\dec}(Z_n)$ resulting from \eqref{eq_solution_star} or \eqref{eq_solution_star_c} will be different in general. A natural task, therefore, is to compare the corresponding regrets of these decisions, i.e., compare $\regret(\optimal{\dec}(Z_n),\theta)$ with $\optimal{\dec}$ defined by either \eqref{eq_solution_star} or \eqref{eq_solution_star_c} for $n=1,\ldots,N$ with some $N\in\nat$.
We refer to \Cref{subsec_example_empmeasure_iid} and \Cref{fig_plot_newsvendor_average_regrets} for such an empirical comparison.

Furthermore, under additional convexity assumptions, which are much more likely to be satisfied by $c$ than by the original $\regret$, the minimum and maximum in the optimal solution can be exchanged (Sion's minimax theorem~\cite{Sion:1958}), which can be advantageous for numerical implementation.
\end{rem}
    
Our main optimization problem~\eqref{eq_optp_main} and its corresponding optimal solution~\eqref{eq_solution_star} can be viewed as a generalization of the distributionally robust optimization framework proposed in~\cite{iP1}, obtained by selecting a specific (degenerate) penalty function~$\ell$. Notably, the choice of penalty function directly influences the structure and properties of the resulting data-driven decision~\eqref{eq_solution_star} such as its consistency. We provide more details in the following remark.

\begin{rem}[Discussion of problem and solution]\label{rem_discussion_problem_solution}
Consider the degenerate loss function $\hat \ell\pcolon [-\infty, \infty]\to[-\infty, \infty]$ defined as $\hat \ell(y) \defeq 0$ if $y\geq 0$ and $\hat \ell(y) \defeq -\infty$ if $y<0$. Notice that in the asymptotic entropic risk measure \Cref{ex_standardLDP_ARE} with $\ell={\hat \ell}$, the constraint in \eqref{eq_optp_main} becomes
\begin{equation}\label{eq_sutter_constraint}
    \limsup_{n\to\infty} \frac{1}{n} \log \probsub{\theta}{\regret(\dec(\proc_n),\theta) > \off(\proc_n)} \leq -\rate,\quad\forall\theta\in\Theta,
\end{equation}
and the respective optimal decision \eqref{eq_solution_star} at $\proc_n$ with the pseudo-inverse of ${\hat \ell}$ given by  $\pseudoinv{{\hat \ell}}(y) = \infty$ if $y\geq 0$ and $\pseudoinv{{\hat \ell}}(y) = 0$ if $y<0$
would under specific conditions (see \cite{iP1}) equal to 
\begin{equation}\label{eq_sutter_solution}
    \argmin_{x\in\dspa} \supb{\theta \in \Theta}{ \regret(x,\theta) - \pseudoinv{{\hat \ell}}\left(\ratefct_\theta(\proc_n) - \rate\right)} = \argmin_{x \in \dspa} \sup_{\theta \in \Theta,\ \ratefct_\theta(\proc_n)< \rate} \regret(x,\theta).
\end{equation}
Hence, the corresponding optimal decision is characterized  by minimizing the regret in $x$, while controlling for the worst choice $\theta\in\Theta$ in a rate-function-induced ``\textit{ball}'' centered around $\proc_n$.
Of course, ${\hat \ell}$ is not an admissible penalty function and we could not prove the optimality or feasibility of \eqref{eq_sutter_solution}, as the problem \eqref{eq_optp_main} would not be well defined in the first place. Nevertheless, the formulation \eqref{eq_sutter_constraint} shows our motivation for $\off$ to be called an upper confidence, as well as for $\rate$ to be called the convergence rate, because heuristically one could rearrange the inequality \eqref{eq_sutter_constraint} and arrive at the interpretation, that the probability of the unfavorable event ``\textit{regret is larger than some as small as possible confidence bound}'' converges to zero exponentially with rate $\rate$. This is exactly in line with our goal of solving \eqref{eq_overall_optp}, i.e., with finding a decision with the smallest possible regret.

Even though this setup is not admissible in our framework, we can easily simulate it. For example, we can choose $\ell=\ell_{\alpha,\beta}$ from \eqref{eq_def_penalty_alphabeta} with large $\beta\gg0$ and small $\alpha\approx0$ to achieve a pointwise approximation of ${\hat \ell}$ (see \Cref{fig_plot_penalty_alphabeta_beta}),
so that the optimization problem with such a fixed penalty function generates a solution approximating \eqref{eq_sutter_constraint}. If we consider the case where $\proc_n$ is an empirical estimator for $\theta$, then \eqref{eq_sutter_solution} is a decision making method
very much related to standard \ac{DRO} approaches, as analyzed by \cite{iP1}, which further underlines the legitimacy of our method. 

\begin{figure}[h!]
    \begin{center}
        \includegraphics[width=0.3\textwidth]{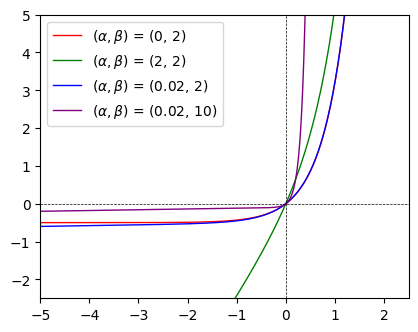}
        \includegraphics[width=0.3\textwidth]{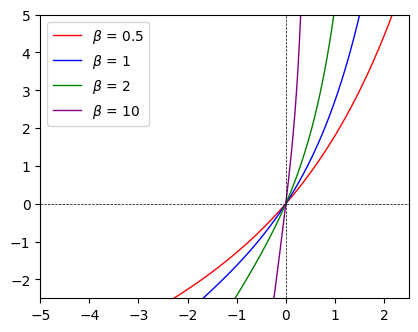}
        \includegraphics[width=0.3\textwidth]{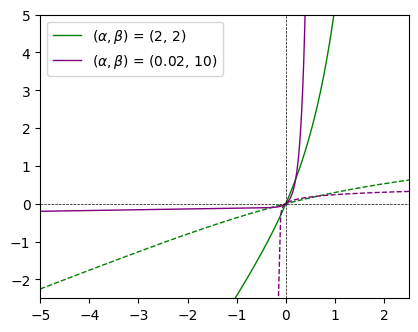}
    \end{center}
    \caption{The plots left and in the middle show inverse penalty functions
    $\inv{\ell}_{\alpha,\beta}(y)= \frac{1}{\beta}(\exp(\beta y) -1) + \alpha y$ and $
    \inv{\ell}_{\beta,\beta}(y)= \frac{1}{\beta}(\exp(\beta y) -1) + \beta y$, $y\in\real$ 
    with $\alpha,\beta>0$. Although $\ell_{\alpha,\beta}$ with $\alpha = 0$ is not an admissible penalty function (since $\inv{\ell}_{0,\beta}>-\frac{1}{\beta}$), a plot of its inverse is still included for illustration.
    The plot on the right displays both the penalty function $\ell_{\alpha,\beta}$ (dashed lines) and its inverse $\inv{\ell}_{\alpha,\beta}$ (solid lines)
    under different choices of $\alpha,\beta>0$.}
    \label{fig_plot_penalty_alphabeta_beta}
\end{figure}

The degenerate case \eqref{eq_sutter_solution} also gives an insight into the consistency properties of the optimal decision in these circumstances, i.e., particularly in situations unlike the one presented in \Cref{ex_plugin_via_LLN}.
Although the consistency of decision methods of type \eqref{eq_sutter_solution} can be achieved in special settings with a decaying radius parameter $\rate_n\to0$, see \cite{ref:Duchi2-21,ref:Mengmeng-21}, it is neither guaranteed nor expected for a constant rate parameter $\rate>0$. This comes as a cost of the robust nature of \eqref{eq_sutter_solution}. Hence, the regret at the evaluated data-driven decision $\optimal{\dec}(\proc_n)$ does not converge to zero as $n\to\infty$, and therefore, it is of interest to study the magnitude of such a \textit{consistency gap}, which we define in the next section, where we also discuss its connection to the penalty function $\ell$.
\end{rem}

\subsection{Consistency}\label{subsec_problem_statement_consistency}
To analyze the consistency properties of the evaluated optimal data-driven decision $\optimal{\dec}(\proc_n)$ in context of the ultimate goal of solving $\eqref{eq_overall_optp}$, we take up on the concept of the \defformat{consistency gap} as already touched upon in the introduction and \Cref{rem_discussion_problem_solution}.
Dealing with such a notion of consistency, it is reasonable to require the process $\seq{\proc}{n}$ to concentrate accordingly to the following assumption, which shall hold throughout this section.
\begin{assumption}[Convergent process]\label{ass_convergent_process} 
There exists a function $\proc_\infty\pcolon \Theta\to\sspa$ such that for every $\theta\in\Theta$ it holds $\lim_{n\to\infty} \proc_n = \proc_\infty(\theta)$ $\mathbb{P}_\theta$-a.s.
\end{assumption}
We define the \defformat{consistency gap} of a decision function $\optimal{\dec}$ as
$\Gamma\pcolon\Theta\to[0,\infty)$ where
\begin{align}\label{eq_def_consistency_gap}
\begin{split}
\theta\mapsto\Gamma(\theta) \defeq \regret(\optimal{\dec}(\proc_\infty(\theta)),\theta) = c(\optimal{\dec}(\proc_\infty(\theta)),\theta) - \min_{y\in\dspa} c(y,\theta),
\end{split}
\end{align}
so that $\Gamma(\theta)=\lim_{n\to\infty} \regret(\optimal{\dec}(\proc_n),\theta)$ holds $\mathbb{P}_\theta$-a.s. for any $\theta\in\Theta$. Further, $\optimal{\dec}$ is called a \defformat{(strongly) consistent} solution to \eqref{eq_optp_main} if and only if $\Gamma(\theta)=0$ for all $\theta\in\Theta$.

Recalling the \Cref{ex_plugin_via_LLN}, it is straightforward to see that the consistency gap of the optimal solution to the \ac{LLN} induced problem~\eqref{eq_optp_main}, namely of the plug-in decision $\decplug$, is always zero.
However, in a general setting with a non-degenerate rate function we cannot expect the decision function $\optimal{\dec}$ to be consistent, and accordingly $\Gamma\neq 0$. Nevertheless, given the \ac{LP}~\eqref{eq_LP} there are other objects besides $\ratefct$ impacting the consistency gap.
With the definition of $\supfct$ in \eqref{eq_G} we obtain
$\supfct(\optimal{\dec}(z),z) \geq \regret(\optimal{\dec}(z),\theta) - \inv{\ell}(\ratefct_{\theta}(z)-\rate)$ for all $\theta\in\Theta$, $z\in\sspa$.
Therefore, the consistency gap at some $\theta\in\Theta$ can be bounded from above by
\begin{align}\label{eq_gap_upperbound}
\begin{split}
    \Gamma(\theta) &\leq \supfct(\optimal{\dec}(\proc_\infty(\theta)),\proc_\infty(\theta)) + \inv{\ell}(\ratefct_{\theta}(\proc_\infty(\theta))-\rate)\\
    &= \min_{x\in\dspa}\max_{\utheta\in\Theta} \left\{ \regret(x,\utheta) -  \inv{\ell}(\ratefct_{\utheta}(\proc_\infty(\theta))-\rate) \right\} + \inv{\ell}(\ratefct_{\theta}(\proc_\infty(\theta))-\rate),
\end{split}
\end{align}
where in the equality we have used the definitions of $\supfct$ and $\optimal{\dec}$.
For typical rate functions, such as the relative entropy $\ent$ in \Cref{ex_main_base} with $\proc_\infty(\theta)=\theta$, it holds $\ratefct_{\theta}(\proc_\infty(\theta))=0$, in which case the last summand in \eqref{eq_gap_upperbound} further simplifies to $\inv{\ell}(-\rate)$.
The condition
\begin{equation}\label{eq_suff_vanishing_gap}
   \min_{x\in\dspa}\max_{\utheta\in\Theta} \left\{ \regret(x,\utheta) -  \inv{\ell}(\ratefct_{\utheta}(\proc_\infty(\theta))-\rate) \right\} + \inv{\ell}(\ratefct_{\theta}(\proc_\infty(\theta))-\rate) = 0,\quad\forall\theta\in\Theta
\end{equation}
is then sufficient for $\Gamma=0$, i.e., for the consistency of $\optimal{\dec}$. However, due to the abstract nature of \eqref{eq_suff_vanishing_gap}, it is not always straightforward to verify it. Nonetheless, as we will discuss in the remainder of this section, the term \eqref{eq_suff_vanishing_gap} highlights that the choice of the penalty function $\ell$ has a major impact on the resulting consistency properties of the optimal and $\ell$-depended decision $\optimal{\dec}$ defined according to \eqref{eq_solution_star}.

In \Cref{rem_discussion_problem_solution} we have outlined how the optimal decision of a problem \eqref{eq_optp_main} with the asymptotic entropic risk measure and penalty function $\ell_{\alpha,\beta}$ with large $\beta\gg0$ and small $\alpha\approx0$ resembles the well known distributionally robust decision method \eqref{eq_sutter_solution}. In such scenarios, the consistency gap $\Gamma$ is usually strictly positive, but the (albeit inconsistent) decision can be still priced for its robustness. Here, the interpretation of the term \textit{robustness} is twofold. 
On the one hand, we refer to the general distributional robustness in the \ac{DRO} sense, which is represented by a ambiguity set defining (pseudo-)metric. In standard \ac{DRO} the choice of such (pseudo-)metric is ambiguous itself, which gets often interpreted as distributional robustness without prior assumptions on $\seq{\proc}{n}$. 
On the other hand, in our case this (pseudo-)metric is given by the rate function (see \eqref{eq_sutter_solution}), and therefore, its choice is tightly connected to the assumptions on the asymptotic behavior of $\seq{\proc}{n}$ characterized by the \ac{LP}~\eqref{eq_LP}. Then, we refer to robustness, which is framed strictly in the set of models $\{\bbP_\theta\mid\theta\in\Theta\}$. This comes in trade off with the now informed choice of the ambiguity set.
In the context of the second interpretation and in the case when $\proc_n$ is an empirical probability measure and $c$ is some expected loss, the \ac{SAA} method (our $\decplug$) stands for the other extreme, thus for straight forward consistency but also for an over-sensitivity to data fluctuations for small sample sizes. It turns out that moving $\alpha$ away from zero can be interpreted as shifting our method on the scale between \ac{DRO} ($\beta\gg0$ \& $\alpha\approx0$) and \ac{SAA} ($\alpha,\beta\gg0$) providing a trade-off between robustness and consistency. In other words, increasing both $\alpha$ and $\beta$ can lead to gradually closing the consistency gap, while still maintaining the flexibility to leave the ambiguity set $\set{\theta \in \Theta}{\ratefct_\theta(\proc_n)< \rate}$ from \eqref{eq_sutter_solution} at the cost of a larger penalty $\inv{\ell}\left(\ratefct_\theta(\proc_n) - \rate\right)$. Such a penalty is the smaller the ``\textit{closer}'' is  $\theta$ to $\proc_n$, so that the maximum in $\optimal{\dec}$ is taken at some $\theta\approx\proc_n$, hence, the similarity to \ac{SAA}.
This dependency between the consistency gap and penalty function helps to understand the influence of the shape of the latter.\footnote{Notice here, that the shared influence of $\ratefct$ and $\ell$ on the consistency properties of $\optimal{\dec}$ is exemplified by the fact that the plug-in decision $\decplug$ can be either exactly revived from our framework via an appropriate \ac{LP} (\Cref{ex_plugin_via_LLN}) or approximated under any \ac{LP} with a well calibrated penalty function (\Cref{subsec_problem_statement_consistency}).}

In summary, supported by the results of our numerical simulations we claim, that the consistency gap $\Gamma$ can be in certain situations shown to converge to zero as $\beta\to\infty$ at any $\theta\in\Theta$.
In the following, we deliver a heuristic argumentation.

Let $\ell_\beta$ be a penalty function parametrized by some $\beta>0$, so that we have the pointwise convergence 
\begin{equation*}
   \lim_{\beta\to\infty} \inv{\ell}_{\beta}(y) = \begin{cases}{}
       +\infty \quad \ & y >0 \\
        0 \quad \ & y=0 \\
        -\infty\quad \ & y<0.
    \end{cases}
\end{equation*}
For example, one might take $\inv{\ell}_{\alpha,\beta}$ with $\alpha=\beta$; however, in sake of simplicity we will proceed with $\inv{\ell}_\beta(y) \defeq \beta y$, $y\in\real$. 
Moreover, assume that $\ratefct_{\utheta}(\proc_\infty(\theta))=0$ if and only if $\theta = \utheta$ and that $\ratefct_{\utheta}(\proc_\infty(\theta))$ grows monotonically as $\theta$ and $\utheta$ are becoming further apart in terms of some metric on $\Theta$. 
For instance in the Gaussian model with a known covariance matrix $\Sigma$ and an unknown mean $\theta$ we have the rate function $\ratefct_{\theta}(z) = \frac{1}{2}(\theta-z)^{\top}\Sigma^{-1}(\theta-z)$ and $\proc_\infty(\theta)=\theta$ for $\theta\in\Theta\subseteq\real^d$, $z\in\sspa=\real^d$ (\cite{iP1}).
Next, we fix $\theta\in\Theta$.
Then, building up on \eqref{eq_gap_upperbound} and recalling that $\min_{y\in\dspa} \regret(y,\theta)=0$ holds and that $\ratefct_{\theta}(\proc_\infty(\theta))=0$ is assumed, we find the following upper bound for the non-negative consistency gap
\begin{align}
\nonumber
   0\leq\Gamma(\theta)&\leq
   \min_{x\in\dspa}\max_{\utheta\in\Theta} \left\{ \regret(x,\utheta) -  \inv{\ell}(\ratefct_{\utheta}(\proc_\infty(\theta))-\rate) \right\} + \inv{\ell}(-\rate)\\ \nonumber
   & = \min_{x\in\dspa}\max_{\utheta\in\Theta} \left\{ \regret(x,\utheta) -  \beta(\ratefct_{\utheta}(\proc_\infty(\theta))-\rate) \right\} - \beta\rate\\ \nonumber
   & = \min_{x\in\dspa}\max_{\utheta\in\Theta} \left\{ \regret(x,\utheta) -  \beta\ratefct_{\utheta}(\proc_\infty(\theta)) \right\} - \min_{y\in\dspa} \regret(y,\theta)\\ \nonumber
   & \leq \sup_{x\in\dspa}\left\{\max_{\utheta\in\Theta} \left\{ \regret(x,\utheta) -  \beta\ratefct_{\utheta}(\proc_\infty(\theta)) \right\} - \regret(x,\theta)\right\}\\
   & = \sup_{x\in\dspa}\left\{\max_{\utheta\in U_\beta(\theta)} \left\{ \regret(x,\utheta) -  \beta\ratefct_{\utheta}(\proc_\infty(\theta)) \right\} - \regret(x,\theta)\right\}, \label{eq_gap_beta_upperbound}
\end{align}
where the second inequality holds by \Cref{lemma_inf_minus_inf_leq_sup} and in the last equality $U_\beta(\theta)\subseteq\Theta$ is a neighborhood around $\theta$, which shrinks as $\beta\to\infty$. Such a sequence of neighborhoods indeed exists, because the penalization term $\beta\ratefct_{\utheta}(\proc_\infty(\theta))$ is strictly positive and increases significantly with $\beta\to\infty$ as long as $\theta\neq\utheta$, and thus, in order to maximize $\utheta\mapsto c(x,\utheta) - \beta\ratefct_{\utheta}(\proc_\infty(\theta))$ one needs to choose ``$\utheta\approx\theta$''. Therefore, as $U_\beta(\theta)$ is getting closer to being the singleton $\{\theta\}$ with increasing $\beta$, we can see that \eqref{eq_gap_beta_upperbound} converges to zero, eventually leading to the consistency gap also approaching zero. 
This becomes obvious with the strictly convex (quadratic) form of the continuous Gaussian rate function as an example. There, we can quantify the neighborhood for any $\beta$ using the continuous derivative of $\ratefct$.
Selecting $\inv{\ell}_\beta=\inv{\ell}_{\beta,\beta}$ from \eqref{eq_def_penalty_alphabeta} displayed in \Cref{fig_plot_penalty_alphabeta_beta} or any similar parametrization instead of the linear penalty function works in analogy.

Note, however, that the optimization problem \eqref{eq_optp_main} is defined with an arbitrary but fixed penalty function and thus, every choice of $\beta$ corresponds to a specific optimization problem with the optimal decision function ``$\optimal{\dec}_\beta$'' as given in \eqref{eq_solution_star}. While a natural approach might be to let $\beta \to \infty$ as $n \to \infty$, this adjustment, which effectively introduces a sequence of penalty functions $\seq{\ell}{n}$, is not supported by our current theoretical framework.
Therefore, in practice one would select a $\beta$ large enough, so that $\Gamma(\theta)\approx 0$ across all $\theta\in\Theta$. 

\section{Generalized Problem Statement}\label{sec_other_opt_criteria}

In this section, we formalize the preceding results and present the most general formulation of our data-driven optimization framework. The key requirement is that the data process must satisfy an abstract Laplace principle. To that end, let $\seq{\proc}{n}$ be a $\sspa$-valued stochastic process on a probability space $(\pspa, \salg, \bbP_\theta)$, where $\sspa$ is a non-empty metric space and $\theta$ is an unknown parameter in a compact, convex metric space $\Theta$. On an intuitive level, as $n \in \mathbb{N}$ increases, the data-generating process reveals progressively more information about the true parameter and should therefore satisfy certain concentration properties.

\subsection{Asymptotic risk measures and concentration} \label{ssec:asymptotic:risk:measures}
Depending on the unknown parameter $\theta \in \Theta$, let $\phi_\theta$ be an asymptotic risk measure mapping continuous transformations $(f(\proc_n))_{n \in \nat}$ of the process $\seq{\proc}{n}$ to the extended real line $\realinf$, and let $\ratefct \pcolon \Theta \times \sspa \to [0, \infty]$ be a rate function meeting \Cref{ass_ratefct}, such that the pair $(\phi_\theta, \mathcal{I}_\theta)$ satisfies the Laplace principle
\begin{equation}\label{eq_LP2}
    \phisub{\theta}{(f(\proc_n))_{n \in \nat}}
    = \sup_{z \in \mathcal{Z}} \left\{ f(z) - \mathcal{I}_\theta(z) \right\},
    \quad \forall f \in \ctsset(\sspa),\ \forall \theta \in \Theta.
\end{equation}
When $\phi_\theta$ represents the standard asymptotic entropic risk measure, the \ac{LP}~\eqref{eq_LP2} constrained to $\ctsset_b(\sspa)$ simplifies to the classical result known as Varadhan's Lemma~\cite[Theorem 4.3.1]{iB1}.
However, this principle is not limited to the asymptotic entropic risk measure $\phi_\theta^{\mathsf{ENT}}$ and there exist other functionals $\psi\pcolon \mathcal{C}(\mathcal{Z}) \to [-\infty, \infty]$ that admit a similar representation. A key property of such functionals is \emph{maxitivity}, meaning that $\psi(f \vee g) = \max\{\psi(f), \psi(g)\}$, where $f \vee g$ denotes the pointwise maximum of the functions $f, g \in \mathcal{C}(\mathcal{Z})$. Maxitive functionals play a central role in possibility theory (\cite{Dubois1988, Zadeh1978}) and in idempotent analysis (\cite{Maslov1997}).

To better illustrate the concept, in the following, we present a list of five typical examples of asymptotic risk measures and their associated concentration properties as characterized by the Laplace principle~\eqref{eq_LP2}, some of which have already been discussed in Section~\ref{sec_problem_statement}. Some of the examples presented below will be revisited in the numerical experiments in Section~\ref{sec_examples}.

\subsubsection{Law of large numbers}\label{subsubsecLLN} We recall Example~\ref{ex_LP_via_LLN}. Suppose that $\Theta = \mathcal{Z}$ and that the sequence $\seq{\proc}{n}$ satisfies an \ac{LLN}, i.e., $\proc_n \to \theta$ as $n \to \infty$ $\bbP_\theta$-a.s.. In that case, we obtain the \ac{LP}~\eqref{eq_LP2} as
\begin{equation*}
    \phiupsub{LLN}{\theta}{(f(\proc_n))_{n \in \nat}} = f(\theta)
    = \sup_{z \in \sspa} \{ f(z) - \ratefct^{\mathsf{LLN}}_\theta(z) \}, \quad f \in \ctsset(\sspa),
\end{equation*}
where the rate function is given by $\ratefct^{\mathsf{LLN}}_\theta(z) = 0$ if $z = \theta$, and $\ratefct^{\mathsf{LLN}}_\theta(z) = \infty$ if $z \neq \theta$. It is straightforward to verify that $\ratefct^{\mathsf{LLN}}$ satisfies Assumption~\ref{ass_ratefct}.

As a concrete example, consider the finite-state i.i.d.~setting introduced in \Cref{ex_main_base}, where the law of large numbers ensures that $Z_n \to \theta$ $\bbP_\theta$-a.s. as $n \to \infty$. Moreover, for an illustration of the resulting optimal decision, see \Cref{ex_plugin_via_LLN}.

\subsubsection{Law of large numbers under distributional uncertainty} \label{ssec:LLN:uncertainty}
The assumption that data at each step are drawn from a fixed distribution is often too strong in practical applications and typically subject to distributional uncertainty. For example, in the finite-state i.i.d.~setting introduced in \Cref{ex_main_base}, the probabilities $\probsub{\theta}{\xi_k = i} = \theta_i$ for $i = 1, \ldots, d$ are often not known precisely. Instead, they are imprecise in the sense of \cite{Walley1991}, meaning they are only specified within some intervals $\probsub{\theta}{\xi_k = i} \in [\underline{\theta}_i, \overline{\theta}_i]$. 
In such cases, the corresponding expectations are not single-valued but instead lie within an interval $[\underline{\mathbb{E}}_\theta,\overline{\mathbb{E}}_\theta]$ that reflects the uncertainty in the underlying model.

In what follows, we focus on the \emph{upper expectation} $\overline{\mathbb{E}}_\theta$, which represents the worst-case expected value over a family of plausible probability distributions (the lower expectation is typically given via the relation $\underline{\mathbb{E}}_\theta[\xi] = -\overline{\mathbb{E}}_\theta[-\xi]$). The upper expectation characterizes the uncertain distribution of $Z_n$ through  $\overline{\mathbb{E}}_\theta[f(Z_n)]$, where $f$ ranges over $\mathcal{C}(\mathcal{Z})$ or a suitable subset on which the upper expectation is well defined.
We assume that, for every $\theta \in \Theta$, the upper expectation $\overline{\mathbb{E}}_\theta$ is defined for all random variables $\xi \pcolon \Omega \to \real$ belonging to a linear space $\mathcal{G}$ such that $c \in \mathcal{G}$ for all constants $c \in \real$ and $|\xi| \in \mathcal{G}$ for all $\xi \in \mathcal{G}$. Moreover, for all $\xi, \zeta \in \mathcal{G}$, we assume:
\begin{enumerate}[label=(\roman*)]
    \item $\overline{\mathbb{E}}_{\theta}[\xi] \leq \overline{\mathbb{E}}_{\theta}[\zeta]$ whenever $\xi \leq \zeta$, \label{item1:sublinear:expec}
    \item $\overline{\mathbb{E}}_\theta[\xi + \zeta] \leq \overline{\mathbb{E}}_\theta[\xi] + \overline{\mathbb{E}}_\theta[\zeta]$, \label{item2:sublinear:expec}
    \item $\overline{\mathbb{E}}_\theta[\xi + c] = \overline{\mathbb{E}}_\theta[\xi] + c$ for all $c \in \real$, \label{item3:sublinear:expec}
    \item $\overline{\mathbb{E}}_\theta[\alpha \xi] = \alpha \overline{\mathbb{E}}_\theta[\xi]$ for all $\alpha \in [0,\infty)$. \label{item4:sublinear:expec}
\end{enumerate}
A functional satisfying these properties is called a \emph{sublinear expectation}, and the triplet $(\Omega, \mathcal{H}, \overline{\mathbb{E}}_\theta)$ is referred to as a \emph{sublinear expectation space} (\cite{iB27}). Sublinear expectations are closely related to coherent risk measures (\cite{Delbaen1999}) and upper coherent previsions in the theory of imprecise probabilities (\cite{Walley1991}). It is well known  that, under suitable continuity assumptions, sublinear expectations admit a robust representation of the form $\overline{\mathbb{E}}_{\theta}[\xi] = \sup_{\mathbb{Q} \in \mathcal{Q}_\theta} \mathbb{E}_{\mathbb{Q}}[\xi]$,
where the ambiguity set $\mathcal{Q}_\theta$ consists of all probability measures regarded as plausible under the parameter $\theta$. If $\mathcal{Q}_\theta$ is a singleton, the model is fully specified; otherwise, it captures distributional uncertainty. Furthermore, we observe that Conditions~\ref{item2:sublinear:expec} and~\ref{item4:sublinear:expec} can be relaxed to a convexity condition, in which case one speaks of a \emph{convex expectation} or a \emph{convex risk measure}. For further details, we refer to \cite{iB28}.

The key assumption is that the sequence $\seq{\proc}{n}$ satisfies a \ac{RLLN}, i.e.,
\[
\lim_{n\to\infty} \overline{\mathbb{E}}_\theta[f(\proc_n)]=\sup_{z\in Z_\infty(\theta)} f(z), \quad \forall f \in \ctsset(\sspa),
\]
for some set $Z_\infty(\theta)\subseteq\mathcal{Z}$. Then, we obtain the \ac{LP}~\eqref{eq_LP2} as
\begin{equation}\label{eq:LP RLLN}
    \phiupsub{RLLN}{\theta}{(f(\proc_n))_{n \in \nat}} \defeq  \lim_{n \to \infty} \overline{\mathbb{E}}_\theta[f(\proc_n)] = \sup_{z \in \sspa} \{ f(z) - \ratefct^{\mathsf{RLLN}}_\theta(z) \}, \quad f \in \ctsset(\sspa),
\end{equation}
where the rate function is given by 
\begin{align}\label{eq:RF RLLN}
   \ratefct^{\mathsf{RLLN}}_\theta(z)  & \defeq \begin{cases}{} 0 & \text{ if }z\in Z_\infty(\theta) \\
   \infty & \text{ otherwise},
    \end{cases} \qquad (\theta,z)\in\Theta\times\sspa.
\end{align}
Whether $\ratefct^{\mathsf{RLLN}}$ satisfies Assumption~\ref{ass_ratefct} depends on the structure of $Z_\infty(\theta)$, which can be verified, for instance, in the following example. For the resulting optimal decision, see \Cref{ex_DRO_via_robustLLN} below. For robust versions of law of large numbers, we refer to \cite{Maccheroni2005, deCooman2008, Peng2008, iB27, Terlan2014, iP6}.

\begin{example}[Finite state i.i.d.~process under distributional uncertainty]\label{ex_main_base:robust}
    Let $\Delta_d$ denote the probability simplex in $\real^d$ for some $d \in \mathbb{N}$, and consider $\Theta = \Delta_d$. 
Let $(\xi_k)_{k \in \mathbb{N}}$ be a finite-state stochastic process, i.e., $\xi_k \pcolon \Omega \to \{1, \ldots, d\}$. For any $\theta \in \Theta$, the law of $\xi_k$ is assumed to lie within a ball of radius $R > 0$ centered at $\theta$, i.e., \[B_\theta(R) \defeq \{ \theta' \in \Theta \mid \mathsf{d}(\theta', \theta) \leq R \},\] where $\mathsf{d}$ is a metric on $\Theta = \Delta_d$. This leads to the intervals $\mathbb{P}_\theta[\xi_k = i] \in \{ \theta'_i \mid \theta' \in B_\theta(R) \}$, $i=1,\ldots,d$. In addition, the random variables $(\xi_k)_{k \in \mathbb{N}}$ are assumed to satisfy a certain independence assumption. Specifically, we require that the distribution of $\xi_k$ is not influenced by conditioning on $\xi_1, \dots, \xi_{k-1}$. That is, we assume that the conditional distribution $\mu_{k-1,k}(\cdot\mid\xi_1, \dots, \xi_{k-1})$ of $\xi_k$ conditioned on $\xi_1, \dots, \xi_{k-1}$ belongs to the set $B_\theta(R)$, and emphasize that $B_\theta(R)$ does not depend on the values of $\xi_1, \dots, \xi_{k-1}$. This means that the sample $\xi_k$ is drawn according to some distribution in $B_\theta(R)$, independently of $\xi_1, \dots, \xi_{k-1}$. In that case, the joint distribution of $(\xi_1,\dots,\xi_k)$ modeled via the respective sublinear expectation is given by
\begin{equation}\label{eq:joint distr}
\overline{\mathbb{E}}_\theta[f(\xi_1,\dots,\xi_k)]=\sup_{\mu}\sum_{(i_1,\dots, i_k)\in \{1,\dots,d\}^k}f(i_1,\dots,i_k)\,\mu(i_1,\dots,i_k),\quad f\pcolon\{1,\dots,d\}^k\to\real,
\end{equation}
where the supremum is taken over all 
$
\mu(i_1, \dots, i_k) = \mu_{0,1}(i_1) \prod_{l=2}^k \mu_{l-1,l}(i_l \mid i_1, \dots, i_{l-1})$ for all conditional distributions \(\mu_{l-1,l} \in B_\theta(R)\). It turns out that the sublinear expectation satisfies
\begin{equation}\label{cond:inpendendence}
\overline{\mathbb{E}}_\theta[f(\xi_1, \dots, \xi_{k})] = \overline{\mathbb{E}}_\theta\left[\overline{\mathbb{E}}_\theta\left[f(i_1, \dots, i_{k-1}, \xi_{k})\right]\big|_{i_1 = \xi_1, \dots, i_{k-1} = \xi_{k-1}}\right]
\end{equation}
for all functions $f\pcolon\{1,\dots,d\}^k\to\real$. This equation reflects that $\xi_{k}$ is independent of $\xi_1, \dots, \xi_{k-1}$. 
Such an independence assumption is, for instance, used in \cite{iB27}.

Finally, let  $\sspa=\Delta_d$ and $\proc_n = \sum_{i=1}^d \left(\frac{1}{n}\sum_{k=1}^{n}\indfctset{\xi_k=i}\right)e_i$ be the empirical measure associated with $(\xi_1,\ldots,\xi_n)$.
As an application of the \ac{LLN} in \cite{Peng2008, iB27} or \cite[Section~3]{iP6}, it follows that for every function $f \pcolon \{1, \dots, d\} \to \real$,
\[
\lim_{n\to\infty}\overline{\mathbb{E}}_\theta[f(Z_n)]=\sup_{z\in B_\theta(R)} f(z).
\]
Hence, $\seq{\proc}{n}$ satisfies the robust law of large numbers with $Z_\infty(\theta)=B_\theta(R)$. Assuming a mild condition on the metric $\mathsf{d}$, the corresponding rate function $\ratefct^{\mathsf{RLLN}}$ satisfies Assumption~\ref{ass_ratefct}.\footnote{Indeed, Assumptions~\ref{ass_ratefct}~\ref{item_ass_ratefct_dom} and \ref{item_ass_ratefct_lsc} always hold. As for \ref{item_ass_ratefct_edgects}, we assume that $\mathsf{d}(\theta,\lambda \theta + (1 - \lambda) \theta')<\mathsf{d}(\theta,\theta')$ for all $\theta\neq\theta'$ and $\lambda\in (0,1]$.
Now, let $(\theta, z) \in \dom_{\Theta \times \sspa} \ratefct^{\mathsf{RLLN}}$ and suppose $z_n \to z$. By definition of the rate function, we have $\mathsf{d}(z, \theta) \le R$, so that, due to the assumption on the metric, there exists a sequence $(\theta_n)_{n \in \mathbb{N}} \subseteq \Theta$ such that $\theta_n \to \theta$ and $\mathsf{d}(z, \theta_n) < R$. Using the triangle inequality of the metric, for all sufficiently large $n$, there exists $\theta_{m_n}$ with $\mathsf{d}(z_n, \theta_{m_n}) \le R$ and $\theta_{m_n} \to \theta$ as $n \to \infty$. This shows that $\ratefct^{\mathsf{RLLN}}_{\theta_{m_n}}(z_n) \equiv 0 \to \ratefct^{\mathsf{RLLN}}_\theta(z) = 0$.}
\end{example}

\subsubsection{Large deviations}\label{subsubsec:LD} 

The best-known setting for a Laplace principle is classical large deviations theory, as explained in detail in Example~\ref{ex_standardLDP_ARE}. Following the notation of \Cref{ex_standardLDP_ARE}, for each $\theta\in\Theta$, let $\phi^{\mathsf{ENT}}_{\theta}$ denote the asymptotic entropic risk measure formally introduced in \eqref{eq_extendedARE}. Under the assumption that the law of $(\proc_n)_{n\in\mathbb{N}}$ satisfies an \ac{LDP}~\eqref{eq_standardLDP} with good rate function $\ent_\theta$, 
the argumentation in Example~\ref{ex_standardLDP_ARE} guarantees the \ac{LP}~\eqref{eq_LP2} as
\begin{equation}\label{eq_extendedARE:new}
    \phiupsub{ENT}{\theta}{(f(\proc_n))_{n\in\nat}}=\sup_{z\in\mathcal{Z}}\{f(z)-\ent_\theta(z)\},\quad \forall f\in\ctsset(\sspa).
\end{equation}

As a concrete example, consider the finite-state i.i.d.~setting introduced in \Cref{ex_main_base}, where the \ac{LP} is ensured by Sanov's theorem. The relative entropy rate function $\ratefct_\theta$ satisfies Assumption~\ref{ass_ratefct}, as shown in \Cref{prop_relentropy_satisfies}. For an illustration of the resulting optimal decision, see \Cref{ex_plugin_via_LLN}. Notably, the setup of \Cref{ex_main_base} serves as an instance of both this and Section~\ref{subsubsecLLN}. Therefore, one can observe how the markedly different shapes of rate functions $\ent$ and $\ratefct^{\mathsf{LLN}}$ resemble the concentration behavior of $\seq{\proc}{n}$ characterized by the respective \ac{LP}, as well as the properties of the respective optimal decisions. The degenerate $\ratefct^{\mathsf{LLN}}$ reflects the statistical consistency of $\decplug$, whereas the more complex $\ent$ mirrors the additional distributional insight reaching beyond \ac{LLN}.

\subsubsection{Large deviations under distributional uncertainty} \label{ssec:LD:uncertainty}
Similar to the law of large numbers discussed in Sections~\ref{subsubsecLLN} and~\ref{ssec:LLN:uncertainty}, large deviations, as introduced in Section~\ref{subsubsec:LD}, can also be analyzed under distributional uncertainty. Following the approach in Section~\ref{ssec:LLN:uncertainty}, we assume that the law of the sequence $(Z_n)_{n\in\mathbb{N}}$ is subject to uncertainty. In the following, we primarily focus on the i.i.d.~setting. There are various ways to model independence under uncertainty; one approach is to assume that the conditional distributions lie within a ball, as outlined in Example~\ref{ex_main_base:robust}. This leads to the independence condition~\eqref{cond:inpendendence}, which has been studied in the context of large deviations, for example, in~\cite{Lacker2020, Eckstein2019, iP6}.

Here, we focus on a stronger notion of independence, as studied in~\cite[Section~5]{iP21}, where a sequence is assumed to be classically independent, but with respect to a family of probability measures. To illustrate this, we consider the setting of Example~\ref{ex_main_base:robust} and focus on the empirical measure $\proc_n = \sum_{i=1}^d \left(\frac{1}{n} \sum_{k=1}^{n} \indfctset{\xi_k = i}\right) e_i$ associated with the sequence $(\xi_k)_{k \in \mathbb{N}}$. As before, the distribution of each $\xi_k\pcolon\Omega\to\{1,\ldots,d\}$ is assumed to lie in a ball $B_\theta(R)$. However, the joint distribution of $(\xi_1, \dots, \xi_k)$ is now given by the sublinear expectation $\overline{\mathbb{E}}_\theta$ in \eqref{eq:joint distr}, where the supremum is taken over all product measures of the form $\mu(i_1, \dots, i_k) = \prod_{l=1}^k \nu(i_l)$ for all $\nu \in B_\theta(R)$. This means that the samples $\xi_k$ are drawn independently according to a single but uncertain distribution $\nu$ from the set $B_\theta(R)$. 

In other words, contrary to the situation in \Cref{ex_main_base} $\seq{\xi}{k}$ is i.i.d. with respect to not just one but any single one element from the set of distributions $B_\theta(R)$. This increases the ambiguity in our assumptions about the underlying data generating process $\seq{\proc}{n}$, which in turn makes our optimal decision method more robust. This is typically useful if the underlying data comes from multiple different but similar sources, as illustrated by our robust newsvendor example in \Cref{subsec_example_empmeasure_iid}.
Another example, in which such situation arises, is when the data are treated as i.i.d., although they are in fact sampled from a time series associated with a problem that is similar, but not identical, to the one being studied, see \cite[Section~6.2]{ref:rychener-2024}. 

In this situation, for suitably chosen parameter space $\Theta$ and the distributional uncertainty quantifying set $B_\theta(R)$, we show\footnote{Alternatively, one could rely on a version of Cram\'er’s theorem \cite[Corollary 5.3]{iP21}; however, for the sake of convenience, we provide a proof for the respective Sanov's version.} in \Cref{prop_robust_finite_sanov} that Sanov's theorem for finite alphabets \cite[Theorem 2.1.10]{iB1} can be reformulated to accommodate such scenarios and thus provide a version of the \ac{LDP} in~\eqref{eq_standardLDP} (see~\eqref{eq_robust_Sanov} for details), with rate function given by the \defformat{robust relative entropy} $\rent \pcolon \Theta \times \sspa \to [0, \infty]$,
\begin{equation}\label{eq_def_robust_relative_entropy}
\rent_\theta(z) \defeq \inf_{\utheta \in B_\theta(R)} \ent_\utheta(z).
\end{equation}
Following the arguments in~\cite[Example~7.1]{iP11}, we obtain that the \defformat{asymptotic robust entropic risk measure} satisfies the \ac{LP}~\eqref{eq_LP2}, as
\[
        \phiupsub{RENT}{\theta}{(f(\proc_n))_{n\in\nat}}\defeq\sup_{m\in\real}\limsup_{n\to\infty}\frac{1}{n}\log\overline{\mathbb{E}}_\theta\left[\exp(n \min\{f(\proc_n),m\})\right]=\sup_{z\in\mathcal{Z}}\{f(z)-\rent_\theta(z)\}
  \]
    for all $f\in\ctsset(\sspa)$. Hence, just as the expectation is replaced by a ``worst-case'' expectation over all measures parameterized by elements of $B_\theta(R)$, the rate function becomes the infimum of the relative entropy over $B_\theta(R)$. To ensure the applicability of our results, the set-valued map $\theta \mapsto B_\theta(R)$ and the parameter space $\Theta$ must be chosen such that the robust relative entropy \eqref{eq_def_robust_relative_entropy} satisfies \Cref{ass_ratefct}. For instance, it suffices to choose any continuous map $\theta \mapsto B_\theta(R)$ such that $B_\theta(R)$ is non-empty, convex and compact, and $\Theta \subseteq \Delta_d \cap (0,1]^d$ (see \Cref{prop_robust_relentropy_satisfies}).

\subsubsection{Maxitive integrals} 
Functionals $\phi_\theta\pcolon \ctsset(\sspa)\to[-\infty,\infty]$ as in~\eqref{eq_LP2}, given by $\phi_\theta(f) \defeq  \sup_{z \in \mathcal{Z}} \left\{ f(z) - \mathcal{I}_\theta(z) \right\}$,
satisfy $\phi_\theta(f) \le \phi_\theta(g)$ whenever $f \le g$, and $\phi_\theta(f + c) = \phi_\theta(f) + c$ for all $c \in \real$. Moreover, as remarked above, they are maxitive and thus admit a representation as a maxitive integral, specifically as a \emph{convex integral} introduced in~\cite{iSP11.2, iSP11.1}.
For instance, it follows from \cite[Theorem 4.2]{iP11}
that
\begin{equation}\label{eq:maxintegral}
\phi_\theta(f) = \sup_{a \in \real} \left\{ a + J_\theta(\{ f > a \}) \right\},\quad\forall f\in \ctsset_{ub}(\mathcal{Z}).
\end{equation}
Here, $J_\theta \pcolon \mathcal{B}(\mathcal{Z}) \to [-\infty, 0]$ is a \defformat{concentration} function, i.e., a set function satisfying $J_\theta(\emptyset) = -\infty$, $J_\theta(\mathcal{Z}) = 0$, and $J_\theta(A) \leq J_\theta(B)$ whenever $A \subseteq B$. Intuitively, $J_\theta$ quantifies the asymptotic concentration of the process $(Z_n)_{n \in \mathbb{N}}$, that is, the plausibility of the process being in a given set. For example, in the case of the law of large numbers (Section~\ref{subsubsecLLN}), where $\proc_n \to \theta$ as $n \to \infty$ $\bbP_\theta$-a.s., we have $J_\theta(A) = 0$ if $\theta \in A$ and $J_\theta(A) = -\infty$ if $\theta \notin A$. Another example of a concentration functional is $J_\theta(A) \defeq \limsup_{n \to \infty} \frac{1}{n} \log \mathbb{P}_\theta[Z_n \in A]$, which captures the exponential rate at which the process $(Z_n)_{n \in \mathbb{N}}$ concentrates in $A$ in the usual sense of large deviations theory. In this case, the corresponding maxitive integral $\phi_\theta$ coincides with the asymptotic entropy (Section~\ref{subsubsec:LD}), and, by the Varadhan–Bryc equivalence (\cite{iB1}), the functional $J_\theta$ satisfies an \ac{LDP} if and only if the maxitive integral $\phi_\theta$ satisfies an \ac{LP}. As shown in \cite{iP10, iP11}, this fundamental equivalence extends to general concentration functionals and their associated maxitive integrals; see, for instance, \cite[Theorem 5.7]{iP11}. For further generalizations in this direction, particularly regarding large deviations theory for non-additive measures and capacities, we refer to \cite{iB29, Nedovic2005, iP23, iP22, iP21, Gao2012} Related contributions can be found in \cite{Eckstein2019, Lacker2020, Backhoff2020}.

Starting with a family of concentrations $J_\theta$, $\theta \in \Theta$, the associated rate function is defined by $\mathcal{I}_\theta(z) \defeq -\inf_{\varepsilon > 0} J_\theta(B_z(\varepsilon))$. If $J_\theta$ satisfies an \ac{LDP}, i.e., $-\inf_{x \in O} \mathcal{I}_\theta(x) \le J_\theta(O)$ for all open sets $O \subseteq \mathcal{Z}$ and $J_\theta(C) \le -\inf_{x \in C} \mathcal{I}_\theta(x)$ for all closed sets $C \subseteq \mathcal{Z}$, which holds if $J_\theta$ satisfies a suitable form of maxitivity and a tightness condition, then the maxitive integral $\phi_\theta$ in \eqref{eq:maxintegral} satisfies the Laplace principle~\eqref{eq_LP2} with rate function $\mathcal{I}_\theta$. 

For a concrete example concerning concentration of capacities, in which the laws of $(Z_n)_{n \in \mathbb{N}}$ are modeled via sublinear expectations, we refer to~\cite[Section~7.1]{iP11}.

\subsection{Generalized problem and solution}\label{subsec_g_and_other}

So far, our article has focused on identifying the optimal decision function $\optimal{\dec}$ and the corresponding upper confidence bound $\optimal{\off}$ according to the optimization problem~\eqref{eq_optp_main}. In this section, we present a generalization of problem~\eqref{eq_optp_main} along with related formulations, and we state the main mathematical results that characterize their optimal solutions. The proofs are postponed to Section~\ref{sec_proofs}. 

As before, for each $\theta \in \Theta$, let $\phi_\theta$ be an asymptotic risk measure, and let $\ratefct \pcolon \Theta \times \sspa \to [0, \infty]$ be a rate function satisfying \Cref{ass_ratefct}. Let $g \pcolon \dspa \times \Theta \to \real$ be a continuous function such that $x \mapsto g(x, \theta)$ is strictly convex for all $\theta \in \Theta$. We focus on the optimization problem $\min_{x \in \mathcal{X}} g(x, \theta)$ and analyze three different data-driven versions of it. Remarkably, all three yield the same optimal solution, indicating a certain stability within the scope of our problem formulations. 
First, we consider our main optimization problem

\begin{equation}\label{eq_optp_main_g} 
\begin{array}{cl}  \displaystyle\inf_{\dec\in\decset, \ \off\in\uset}& \{\off( z)\}_{ z\in\sspa}\\  
    \displaystyle \ \subjectto &\phisub{\theta}{\left(\ell\left(g(\dec(\proc_n),\theta)-\off(\proc_n)\right)\right)_{n\in\nat}}\leq -\rate,\quad \forall \theta\in\Theta,
\end{array} 
\end{equation}
where $\rate>0$ is the convergence rate parameter and $\ell\pcolon \real\map\real$ is a
penalty function, i.e., an increasing continuous and bijective function with $\ell(0)=0$, which has a continuous extension  $\ell\pcolon [-\infty,\infty]\map[-\infty,\infty]$ by setting $\ell(\infty) \defeq \infty$ and $\ell(-\infty) \defeq -\infty$.
We define the generalized counterparts to \eqref{eq_solution_star} as
\begin{align} \label{eq_solution_star_g}
\begin{split}
    \optimal{\dec}(z) &\defeq \argmin_{x\in\dspa} \supfct(x,z) \quad\text{and}\\
    \optimal{\off}(z) &\defeq \min_{x\in\dspa}\supfct(x,z) = \supfct(\optimal{\dec}(z),z), \quad z \in \sspa,
\end{split}
\end{align}
where 
\begin{align}\label{eq_G_g}
  \supfct(x,z)\defeq \maxb{\theta \in \Theta}{ g(x,\theta) - \inv{\ell}\left(\ratefct_\theta(z) - \rate\right)}, \quad z \in \sspa, x \in \dspa.
\end{align}
It turns out that the functions $\optimal{\dec}$, $\optimal{\off}$, and $\supfct$ are continuous, and that the pair $(\optimal{\dec}, \optimal{\off})$ solves the optimization problem~\eqref{eq_optp_main_g}. This is formalized in the following theorem.

\begin{thm}[Feasibility and optimality of $(\optimal{\dec},\optimal{\off})$]\label{thm_feas_dom_main_g} 
    Let the pair $(\phi_\theta,\ratefct_\theta)$ satisfy the \ac{LP}~\eqref{eq_LP2} and let \Cref{ass_ratefct} hold. Then, the following is true.
    \begin{enumerate}[label=\enumthm]
        \item \label{item_thm_feas_dom_main_feas_g} The pair $(\optimal{\dec}, \optimal{\off})$ in \eqref{eq_solution_star_g} is a feasible solution to \eqref{eq_optp_main_g}, i.e., $(\optimal{\dec},\optimal{\off})\in\decset\times\uset$ and 
        \begin{equation*}
            \phisub{\theta}{\left(\ell\left(g(\optimal{\dec}(\proc_n),\theta)-\optimal{\off}(\proc_n)\right)\right)_{n\in\nat}}\leq -\rate,\quad \forall \theta\in\Theta.
        \end{equation*}
        \item \label{item_thm_feas_dom_main_dom_g} The pair $(\optimal{\dec},\optimal{\off})$ is \defformat{optimal} in \eqref{eq_optp_main_g}, i.e., for any other feasible pair $(\dec, \off)$, it holds $$\optimal{\off}(z)\leq\off(z),\quad \forall z\in\sspa.$$
    \end{enumerate}
\end{thm}

In the special case where $g = \regret$, the optimization problem~\eqref{eq_optp_main_g} reduces to~\eqref{eq_optp_main}, so~\Cref{thm_feas_dom_main_g} includes \Cref{thm_feas_dom_main} as a special case. In an analogy, this also holds for the choice $g=c$, as discussed in \Cref{rem_g_is_c}.

\begin{example}[DRO decision as optimal solution in RLLN context] \label{ex_DRO_via_robustLLN}
In analogy to \Cref{ex_plugin_via_LLN}, we briefly discuss the optimal solution of \Cref{thm_feas_dom_main_g} in the context of 
Section~\ref{ssec:LLN:uncertainty}, which addresses the law of large numbers under distributional uncertainty. This setting naturally lends itself to an interpretation of distributionally robust decision making. Due to the presence of distributional uncertainty, the stochastic process $(\proc_n)_{n\in\mathbb{N}}$ is, roughly speaking, assumed to converge asymptotically to a set $Z_\infty(\theta) \subseteq \Theta$, rather than to a single point in $\Theta$. Suppose that the pair $(\phi_\theta^{\mathsf{RLLN}}, \ratefct^{\mathsf{RLLN}}_\theta)$, as defined in \eqref{eq:LP RLLN} and \eqref{eq:RF RLLN}, satisfy the \ac{RLLN} based \ac{LP}~\eqref{eq:LP RLLN}. In this setting, the optimal decision problem \eqref{eq_solution_star_g} simplifies to
\begin{align}\label{ex:DRO:general}
    \dec^\star(z) = \argmin_{x\in\dspa} \max\{ g(x,\theta)  \mid \theta\in\Theta, \  z\in Z_\infty(\theta) \},\quad z\in\sspa.
\end{align}
To illustrate the connection of the decision~\eqref{ex:DRO:general} to \ac{DRO} more concretely, we consider the finite-state i.i.d.~setting described in \Cref{ex_main_base:robust} with $g=c$. There, the distributional uncertainty is modeled via a ball of radius $R > 0$ centered at $\theta$, i.e., $Z_\infty(\theta) = B_\theta(R)= \{\theta'\in \Theta \ : \ \mathsf{d}(\theta',\theta)\leq R\}$, where $\mathsf{d}$ is a metric on $\Theta$ that is conform with \Cref{ass_ratefct}~\ref{item_ass_ratefct_edgects}. In this case, the optimal decision~\eqref{ex:DRO:general} reduces to the classical form of a \ac{DRO} problem
\begin{equation*}
    \dec^\star(z) = \argmin_{x\in\dspa} \max \{ c(x,\theta)  \mid  \theta\in B_z(R) \},\quad z\in\sspa.
\end{equation*}
Thus, if we work with an \ac{LP}, which is designed to fully (and only) represent either the \ac{LLN} or \ac{RLLN} assumption on the underlying data generating process, then \Cref{thm_feas_dom_main_g} certifies either the plug-in (see \Cref{ex_plugin_via_LLN}) or \ac{DRO} method as optimal decision respectively. The added uncertainty represented by the shift from \ac{LLN} to \ac{RLLN} is also mirrored in the shift from consistency to robustness of the resulting decision function.
\end{example}

\begin{rem}[Bounds on rate function] \label{rem_bounds} 
    In case that the rate function $\ratefct$ in \eqref{eq_LP2} is not known analytically, but can be estimated from below by a rate function $\tilde{\ratefct}$ satisfying \Cref{ass_ratefct}, the candidate pair $(\optimal{\dec}, \optimal{\off})$ from~\eqref{eq_solution_star_g} remains a feasible solution to~\eqref{eq_optp_main_g}. In other words, the upper confidence bound $\optimal{\off}$ still satisfies the guarantees imposed by the feasibility constraints. However, $\optimal{\off}$ is not necessarily optimal in the sense that there may exist a smaller upper confidence bound that also satisfies the feasibility constraint. This follows directly from the proof of \Cref{thm_feas_dom_main_g}-\ref{item_thm_feas_dom_main_feas}, where the \ac{LP}~\eqref{eq_LP2} is replaced by the inequality 
    $\phisub{\theta}{(f(\proc_n))_{n \in \nat}} \leq \sup_{z \in \sspa} \{ f(z) - \tilde{\ratefct}_\theta(z) \}$ for all $f\in\ctsset(\sspa)$, $\theta\in\Theta$.
\end{rem}

\begin{rem}[Consistency]\label{rem_consistency_g}
    The content of \Cref{subsec_problem_statement_consistency} applies in the context of the optimization problem \eqref{eq_optp_main_g} as well. That is, if we generalize~\eqref{eq_def_consistency_gap} to $\Gamma\pcolon\Theta\to[0,\infty)$, $$\Gamma(\theta) \defeq g(\optimal{\dec}(\proc_\infty(\theta)),\theta) - \min_{y\in\dspa} g(y,\theta)$$ while still complying with \Cref{ass_convergent_process}, the (heuristic) analysis of consistency properties in \Cref{subsec_problem_statement_consistency}, which effectively represents the special cases $g=\regret$ and $g=c$, can be done with any other $g$ in a complete analogy. Then, $\optimal{\dec}$ defined in \eqref{eq_solution_star_g} is called a \defformat{consistent} solution to \eqref{eq_optp_main_g} if and only if $\Gamma(\theta)=0$ for all $\theta\in\Theta$ and the sufficient condition for consistency \eqref{eq_suff_vanishing_gap} translates to
    \begin{equation*}
        \min_{x\in\dspa}\max_{\utheta\in\Theta} \left\{ g(x,\utheta) -  \inv{\ell}(\ratefct_{\utheta}(\proc_\infty(\theta))-\rate) \right\} + \inv{\ell}(\ratefct_{\theta}(\proc_\infty(\theta))-\rate) - \min_{y\in\dspa}g(y,\theta)= 0,\quad\forall\theta\in\Theta,
    \end{equation*}
    where $\ratefct_{\theta}(\proc_\infty(\theta))$ is typically zero.
\end{rem}

A closer look at the optimizers~\eqref{eq_solution_star_g} of~\eqref{eq_G_g} reveals that the optimization problem~\eqref{eq_optp_main_g} is closely related to the following reformulation:
\begin{align}\label{eq_optp_modified_argmin_g}
    \begin{array}{cl} \displaystyle 
    \inf_{\modsupfct\in\tilde{\ctsset}(\dspa\times\sspa)}
    &\{\displaystyle\min_{x\in\dspa}\modsupfct(x, z)\}_{ z\in\sspa}\\  
        \ \subjectto &\
      \displaystyle  
      \phisub{\theta}{\left(\ell\Big( g(\dec(\proc_n),\theta)-\min_{x\in\dspa}\modsupfct(x,\proc_n)\Big)\right)_{n\in\nat}}\leq -\rate,\quad \forall \theta\in\Theta,\\
       & \ \text{with } \dec(\cdot)\defeq\argmin_{x\in\dspa} \modsupfct(x,\cdot),
   \end{array}
\end{align}
where $\tilde{\ctsset}(\dspa\times\sspa)\subseteq\costset$ is the set of all jointly continuous functions $\modsupfct\pcolon \dspa\times\sspa\to\real$, so that $x\mapsto \modsupfct(x,z)$ is strictly convex for all $z\in\sspa$.
Since $G$ is strictly convex in the first argument, we obtain from \Cref{thm_feas_dom_main_g} that the optimal solution of \eqref{eq_optp_modified_argmin_g} is given by $\optimal{H}\defeq G$ with optimal decision function  $\optimal{\dec}(\cdot)\defeq\argmin_{x\in\dspa} \supfct(x,\cdot)$
and optimal upper confidence bound $\optimal{\off}(\cdot)\defeq \min_{x\in\dspa} \supfct(x,\cdot)$.
This shows that the optimization problems \eqref{eq_optp_main_g} and \eqref{eq_optp_modified_argmin_g} are strongly linked.
For example, in the case $g=c$ this problem characterizes functions $\modsupfct$ that induce decision functions $\dec(z)=\argmin_{x\in\dspa} \modsupfct(x,z)$, $z\in\sspa$, such that the behavior of the corresponding out-of-sample cost $c(\dec(\proc_n), \theta)$ is controlled. Intuitively, $\min_{x\in\dspa}\modsupfct(x,\proc_n)$ is a minimal over-approximation of  $c(\dec(\proc_n), \theta)$.
In this regard, problem~\eqref{eq_optp_modified_argmin_g} is inspired by~\cite{iP1}.

Furthermore, we present a third formulation, which further highlights the robustness inherent in the structure of our solution. The corresponding optimization problem is
\begin{align}   \label{eq_optp_modified_infinf_g}
    \begin{array}{cl}  \displaystyle \inf_{\modsupfct\in\costset} &\{\displaystyle\min_{x\in\dspa}\modsupfct(x, z)\}_{ z\in\sspa}\\ 
       \ \subjectto & 
       \displaystyle
      \phisub{\theta}{\left(\ell\Big(\min_{x\in\dspa} g(x,\theta)-\min_{x\in\dspa}\modsupfct(x,\proc_n)\Big)\right)_{n\in\nat}}\leq -\rate,\quad \forall \theta\in\Theta.
   \end{array}
\end{align}
If we were to choose $g=\regret$, this optimization problem, unlike problem~\eqref{eq_optp_main}, would aim to find a function $\modsupfct$ such that $\min_{x \in \dspa} \modsupfct(x, Z_n)$ is the least conservative upper confidence bound to the optimal regret $\min_{x\in\dspa} \regret(x,\theta)=0$, so that in this problem formulation the upper confidence bound $\off$ obtains a more specific structure.  Here, ``least conservative'' is expressed via the objective function, and the notion of the upper confidence bound refers to the constraint. 
In contrast to \eqref{eq_optp_main_g}, where the goal is to directly determine an optimal decision function \( \dec \) and characterize the behavior of \( g(\dec(\proc_n),\theta) \) through its constraint, the problem \eqref{eq_optp_modified_infinf_g} should be interpreted as providing the least conservative upper confidence bound on $\min_{x\in\dspa}g(x,\theta)$. This upper confidence property is conveyed through the constraint but does not directly yield a data-driven decision. Under additional convexity assumptions, we show that $\supfct$ 
is an optimal solution in \eqref{eq_optp_modified_infinf_g} as well.

\begin{proposition}[Solution to \eqref{eq_optp_modified_infinf_g}]\label{prop_feas_dom_modified_infinf_g}
    Let the pair $(\phi_\theta,\ratefct_\theta)$ satisfy the \ac{LP}~\eqref{eq_LP2} and let \Cref{ass_ratefct} hold. Suppose that the function $\theta\mapsto g(x,\theta) - \inv{\ell}(\ratefct_\theta ( z)-\rate)$ is concave on $\Theta$ for all $x\in\dspa$, $z\in\sspa$.
    Then, the following holds.
    \begin{enumerate}[label=\enumproposition]
        \item $\supfct$ is feasible in \eqref{eq_optp_modified_infinf_g}, i.e., $\supfct\in\costset$ and 
        \begin{equation*}
            \phisub{\theta}{\left(\ell\Big(\min_{x\in\dspa} g(x,\theta)-\min_{x\in\dspa}\supfct(x,\proc_n)\Big)\right)_{n\in\nat}}\leq -\rate,\quad \forall \theta\in\Theta.
        \end{equation*}
        \item \label{item_prop_feas_dom_modified_infinf_dom} $\supfct$ is optimal in \eqref{eq_optp_modified_infinf_g}, i.e., 
        for any other feasible candidate $\modsupfct$ it holds 
        \begin{equation*}
            \min_{x\in\dspa} \supfct(x,z) \leq \min_{x\in\dspa} \modsupfct(x,z), \quad \forall z\in\sspa.
        \end{equation*}
    \end{enumerate}
\end{proposition}

\begin{rem}[Saddle point]\label{rem_concave}
The concavity assumption in Proposition~\ref{prop_feas_dom_modified_infinf_g}, together with our standing assumptions, ensures that for every $z \in \mathcal{Z}$, the function $\mathcal{X} \times \Theta \to [-\infty, \infty)$ given by
\begin{equation}\label{eq:saddlepoint}
    (x, \theta) \mapsto g(x, \theta) - \inv{\ell}(\ratefct_\theta(z) - \rate)
\end{equation}
is convex and \ac{lsc} in the first argument, and concave and \ac{usc} in the second argument. The latter follows from the proof of Lemma~\ref{lemma_optimal_G_g_cts}. Hence, by Sion's minimax theorem~\cite{Sion:1958}, the function in~\eqref{eq:saddlepoint} admits a saddle point $(\optimal{x}, \optimal{\theta})$. Although the function in~\eqref{eq:saddlepoint} is not real-valued in general, it is possible to restrict the second argument to a compact subset $\tilde\Theta \defeq  \{\theta \in \Theta \mid \mathcal{I}_\theta(z) \le m\}$ for some sufficiently large $m \in \real$, since $g(x, \theta)$ is bounded over $\mathcal{X} \times \Theta$. This restriction makes the function real-valued and preserves the same saddle point.

 For instance, the concavity assumption holds in the case $g = c$, where $c$ is a cost function that is concave in $\theta$, provided that the rate function $\mathcal{I}$ is convex in $\theta$. In contrast, the choice $g = \regret$ is typically not concave in $\theta$.
\end{rem}
\section{Modeling Examples and Numerical Experiments}\label{sec_examples}
This section introduces concrete modeling examples
and numerically examines two classical problem classes from Operations Research--the newsvendor problem (\Cref{subsec_example_empmeasure_iid}) and an optimal portfolio selection problem (Section~\ref{subsec_example_empmean_iid}). In order to ensure the applicability of \Cref{thm_feas_dom_main} and \Cref{thm_feas_dom_main_g} in all these situations, we especially confirm the validity of \Cref{ass_ratefct} while delegating detailed proofs to the \nameref{sec_appendix}.
These examples serve to empirically demonstrate the performance of our method and more importantly to illustrate how our general framework, grounded in the concept of sublinear expectations, can be applied to stochastic optimization problems under additional uncertainty, a setting that is not addressed by existing related approaches such as \cite{iP1}.
Furthermore, we provide empirical evidence highlighting the influence of the shape of the penalty function $\ell$ on the consistency gap $\Gamma$, as introduced in Section~\ref{subsec_problem_statement_consistency} and \Cref{rem_consistency_g}.

Throughout this section, we consider a penalty function, characterized by its inverse
\begin{equation*}
    \inv{\ell}_{\beta}(y)\defeq \frac{1}{\beta}(\exp(\beta y) -1) + \beta y ,\quad y\in\real,
\end{equation*}
with some $\beta>0$, i.e., a special case of \eqref{eq_def_penalty_alphabeta}. The asymptotic risk measure $\phi_\theta$ considered is either the asymptotic entropic risk measure $\phi_\theta^{\mathsf{ENT}}$ defined in \eqref{eq_extendedARE} or its robust counterpart as introduced in \Cref{ssec:LD:uncertainty}. 
Recall, the optimal decision function $\optimal{\dec}$ is given by  \eqref{eq_solution_star_g} as 
\begin{equation*}
    \optimal{\dec}(z) = \argmin_{x\in\dspa} \maxb{\theta \in \Theta}{ g(x,\theta) - \inv{\ell}_\beta\left(\ratefct_\theta(z) - \rate\right)}, \quad z\in\sspa.
\end{equation*}
so that, any selection of the pair $(\phi_\theta,\ratefct_\theta)$, $g\in\{\regret,c\}$ and the parameters $\beta,\rate>0$ defines a specific decision, which is optimal in \eqref{eq_optp_main_g}. The same applies to the definition of the consistency gap $\Gamma$ of such a decision in the spirit of \Cref{rem_consistency_g}.
Thus, an empirical comparison across multiple such selections is desirable.

\subsection{Empirical measure of finite state i.i.d.~process}\label{subsec_example_empmeasure_iid}
The first modeling class we consider is that of a finite-state i.i.d.~process, as introduced in \Cref{ex_main_base}. To demonstrate the versatility of the framework developed in this paper, we examine two widely used asymptotic risk measures and their associated concentration behavior: the asymptotic entropic risk measure (see \Cref{sssec:as:entropic:risk}) and its robust counterpart (see \Cref{sssec:rob:asym:ent:risk}).\footnote{Code for the experiments presented in this section is available at \url{https://github.com/radeksalac/Asymptotic-Optimality-in-Data-Driven-Decision-Making}}
For both these asymptotic risk measures, our focus is on a Laplace principle of the form \eqref{eq_LP2}, derived via an \ac{LDP}. In other words, we consider the settings described in \Cref{subsubsec:LD} and \Cref{ssec:LD:uncertainty}.

\subsubsection{Asymptotic entropic risk measure}\label{sssec:as:entropic:risk}
Consider the setting of \Cref{ex_main_base}
and let $\seq{\xi}{k}$ be an i.i.d.~stochastic process with the finite state space $\Xi\defeq\{1,\ldots,d\}$ and $\probsub{\theta}{\xi_k = i}= \theta_i$ for $i=1,\ldots,d$, where $\theta\in\Theta\subseteq\Delta_d$ is an unknown vector from the probability simplex. Moreover, let $\seq{\proc}{n}$ be the empirical measure process with $\sspa=\Delta_d$, i.e., $\proc_n = \sum_{i=1}^d \left(\frac{1}{n}\sum_{k=1}^{n}\indfctset{\xi_k=i}\right)e_i$. Then, as established in \Cref{ex_main_base}, the asymptotic entropic risk measure $\phi^{\mathsf{ENT}}_{\theta}$ constitutes the \ac{LP}~\eqref{eq_LP2} together with the relative entropy rate function $\ent_\theta$ defined in \eqref{eq_def_relative_entropy}.
That is, we are in the setting of \Cref{subsubsec:LD}.
Indeed, we can show that our main results \Cref{thm_feas_dom_main} and \Cref{thm_feas_dom_main_g} with appropriate cost function and parameter space $\Theta$ hold, as \Cref{ass_ratefct} is satisfied.
\begin{proposition}\label{prop_relentropy_satisfies}
    The relative entropy $\ent$ defined in \eqref{eq_def_relative_entropy} satisfies \Cref{ass_ratefct}
    with $\sspa=\Delta_d$ and either $\Theta=\Delta_d$ or some convex compact $\Theta\subseteq\Delta_d\cap(0,1]^d$; in fact, in the latter case is $\ent$ even continuous.
\end{proposition}
A proof of \Cref{prop_relentropy_satisfies} is provided in the \nameref{sec_appendix}.
Moreover, by the strong law of large numbers it holds $\proc_{n} \to \theta$ $\bbP_\theta$-a.s. as $n\to\infty$, so that \Cref{ass_convergent_process} is satisfied with $\proc_\infty(\theta) = \theta$, $\theta\in\Theta$.
In the following, we demonstrate the performance of the proposed data-driven decision in a well-known setting from Operations Research--the newsvendor problem.

\paragraph{Newsvendor problem.} 
We examine a seller of a perishable good, which loses all value at the end of each day. At the start of every day, the newsvendor places an order of $x\in\dspa$ units from a supplier, paying a wholesale price $\kappa\geq 0$ per unit, where $\dspa=[0,d]$ represents the possible order quantities, here we assume that it is possible for the newsvendor to order a continuous amount of the good. 
Also we assume that the newsvendor faces a transportation cost that is quadratic in the amount of orders he carries, expressed as $\rho x^2$, where $\rho>0$ is some transportation cost parameter.
Afterward, the uncertain discrete demand $\xi$ with state space $\{1,\dots,d\}$ is realized according to an unknown distribution $\theta\in\Delta_d$. The newsvendor sells the good at a retail price $p>\kappa$ per item, continuing until either the demand is met or the inventory is depleted. The newsvendor problem can be formulated as an instance of a finite state i.i.d.~process with expected loss.
The loss (negative reward) of the newsvendor therefore, can be formulated as $h(x,\xi) = \kappa x + \rho x^2 - p\min\{x,\xi\}$ leading to a total expected cost 
\begin{equation}\label{eq_newsvendor_totalexpectedcost}
    c(x,\theta) = \expecsub{\theta}{h(x,\xi)} = \sum_{i=1}^d h(x,i)\theta_i,\quad (x,\theta)\in\dspa\times\Theta.
\end{equation}
Clearly, the cost function $c$ is continuous and strictly convex in $x$ for all $\theta\in\Theta$. 
Altogether, this scenario fits the setting described above and both \Cref{thm_feas_dom_main} as well as \Cref{thm_feas_dom_main_g} are applicable.

\begin{figure}[h!]
    \begin{center}
        \includegraphics[width=1\textwidth]{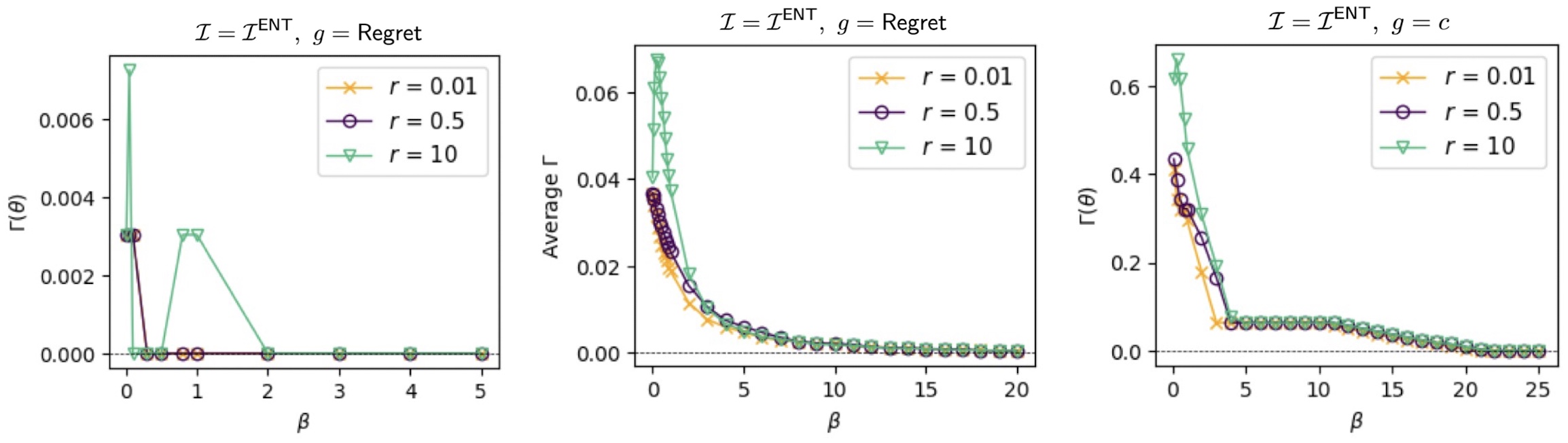}
    \end{center}
    \caption{Plot of the (average) consistency gap of $\optimal{\dec}$ as defined by \eqref{eq_solution_star_g} with $\ratefct=\ent$, $g\in\{\regret,c\}$ and $\rate\in\{0.01,0.5,10\}$ viewed as a function of the penalty function parameter $\beta$ within the standard Newsvendor problem. In the left and right plots, we display consistency gap at the true $\theta$, i.e., $\Gamma(\theta)=\regret(\optimal{\dec}(\proc_\infty(\theta)),\theta)$, whereas in the middle one, we show an average gap across multiple unknown $\utheta_k$'s, i.e., $\frac{1}{200}\sum_{k=1}^{200}\Gamma(\utheta_k)$, and the $\utheta_k$'s are sampled uniformly from $\Delta_d\cap[0.05,1]^d$. The total expected cost \eqref{eq_newsvendor_totalexpectedcost} with $k = 1$, $p = 1.65$, $\rho = 0.0025$ and data generating process parameters $d = 8$, $\Theta = \Delta_8\cap[0.001,1]^8, $ $\theta = (0.115, 0.115, 0.115, 0.125, 0.135, 0.135, 0.135, 0.125)^\top$ are considered.
    }
    \label{fig_plot_newsvendor_GapBeta}
\end{figure}

We consider a synthetic example with an in practice unknown distribution $\theta$ on $\Xi$. In \Cref{fig_plot_newsvendor_GapBeta}, we show the consistency gap of $\optimal{\dec}$ at this true $\theta$, i.e., $\Gamma(\theta)=\regret(\optimal{\dec}(\proc_\infty(\theta)),\theta)$, plotted against the penalty function parameter $\beta$. We do this for both $g=\regret$ and $g=c$ and for different values of $\rate$. In other words, we can observe the performance of the decision $\optimal{\dec}$ under different parameter selections after the law of large number satisfying process $\seq{\proc}{n}$ has converged.
In analogy, we also show the average gap across multiple unknown $\utheta_k$'s, i.e., $\frac{1}{200}\sum_{k=1}^{200}\Gamma(\utheta_k)$, where $\utheta_k$'s are sampled uniformly from $\Delta_d\cap[\epsilon,1]^d$.

 Aligned with the heuristics in \Cref{subsec_problem_statement_consistency}, we empirically confirm that the consistency gap $\Gamma(\theta)$ at the unknown $\theta$ saturates at zero once $\beta$ becomes large enough (\Cref{fig_plot_newsvendor_GapBeta}). Naturally, the magnitude of such sufficiently large $\beta$ stands in relation to $\sup_{x}\abs{g(x,\theta)}$, so that the case $g=\regret$ likely requires smaller $\beta$ than $g=c$ to define an optimization problem \eqref{eq_optp_main_g} with a consistent optimal decision $\optimal{\dec}$, i.e., $\optimal{\dec}$ with $\Gamma(\theta)\approx0$. Notably, $\sup_{x}\abs{g(x,\theta)}$ changes with $\theta$ and so does the point of saturation. We can see this while comparing the plot of consistency gap with $g=\regret$ at the true $\theta$, where the gap is approximately zero for $\beta$ as small as 0.3, and the corresponding plot of averaged consistency gap, where a much larger $\beta$ is needed to achieve similar results. Moreover, our experiments also showed, that the averaged consistency gap decreases as the lower bound $\epsilon$ for $\probsub{\utheta}{\xi_k=i}$ for any $i\in\Xi$ increases. Finally, it does not come as a surprise, that a larger rate $\rate$ and thereby a stronger penalization by the term $\inv{\ell}_\beta(\ent_\theta(z)-\rate)$ in the formula for $\optimal{\dec}$ leads to a larger consistency gap for smaller $\beta$. For $\beta\gg0$ this effect is diminished, as it becomes irrelevant in comparison to the more and more rapidly increasing penalty function $\ell_\beta$.

\begin{figure}[h!]
    \begin{center}
        \includegraphics[width=1\textwidth]{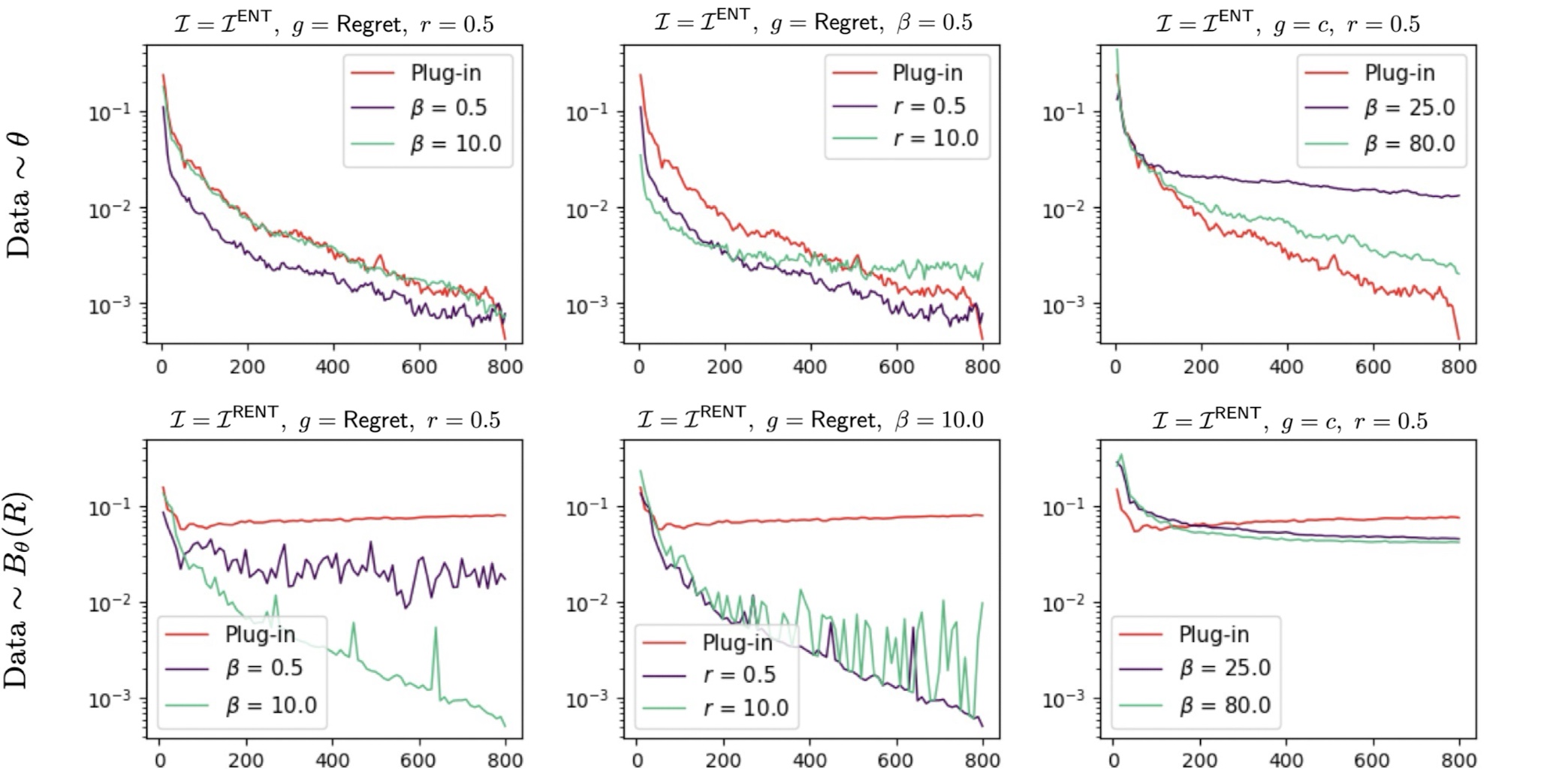}
    \end{center}
    \caption{
    The plots show the averaged regrets for various data-driven decision functions $\dec$ within the (robust) Newsvendor problem. That is, for each sample size $n\leq N$ on the x-axis there is an average regret value $\frac{1}{K}\sum_{k=1}^K\regret(\dec(z^{(k)}_{n}),\theta)$ on the log-scaled y-axis. The results for the plug-in method, i.e., $\dec=\decplug$, are displayed in red, so that the colors green and violet correspond to our decision method, i.e.,  $\dec=\optimal{\dec}$ as defined by \eqref{eq_solution_star_g} with $\ratefct\in\{\ent,\rent\}$, $g\in\{\regret,c\}$ and some selection of parameters $\beta,\rate>0$. The total expected cost \eqref{eq_newsvendor_totalexpectedcost} with $k = 1$, $p = 1.65$, $\rho = 0.0025$ and the data generating process parameters $d = 8$, $\Theta = \Delta_8\cap[0.001,1]^8$,
    $\theta = (0.115, 0.115, 0.115, 0.125, 0.135, 0.135, 0.135, 0.125)^\top$, $N = 800$ and $K = 300$ are considered. 
    The first row of plots belongs to the standard Newsvendor problem, and hence, the data is generated according to this unique and in practice unknown $\theta$. To generate data for the robust Newsvendor in the second row, we consider the closed $L_2$-ball intersected with $\Theta$, $B_\theta(R)$, with $R = 0.07$ and $\utheta_1 = (0.109, 0.145, 0.155, 0.135, 0.12, 0.12, 0.10, 0.116)^\top$, $\utheta_2 = (0.088, 0.09, 0.09, 0.105, 0.152, 0.16, 0.16, 0.155)^\top$, $\utheta_1,\utheta_2\in B_\theta(R)$, with the sampling probability $p_{\utheta_1} = 0.5$.
    }
    \label{fig_plot_newsvendor_average_regrets}
\end{figure}

Next we generate $K$ i.i.d.~samples of length $N$ from a distribution $\theta$ on $\Xi$, i.e., $K$ realizations of $(\xi_1,\ldots,\xi_N)$. This then translates to $K$ sample paths of the empirical distribution process $\seq{\proc}{n}$ of the same length $N$, denoted as $(z_1^{(k)},\ldots,z_N^{(k)})$ with $k=1,\ldots,K$. In the first row of \Cref{fig_plot_newsvendor_average_regrets}, the average regret of the decision $\optimal{\dec}$ under the unknown true distribution $\theta$ based on a sample of sizes $n=5,10,15,\dots,N$ across $K$ experiments is plotted for $g\in\{c,\regret\}$ and for different values of $\beta$ and $\rate$. The simple plug-in decision \eqref{eq_plug} is included as a reference. 

Also here, the general behavior supports our heuristics and confirms our expectations as outlined in \Cref{rem_discussion_problem_solution} and \Cref{subsec_problem_statement_consistency}. Namely, increasing $\beta$ enhances the consistency properties of our decision $\optimal{\dec}$, which thereby produces results similar to $\decplug$. Furthermore, the optimal decision $\optimal{\dec}$ with $g=\regret$ outperforms not just $\optimal{\dec}$ with $g=c$, but also the plug-in decision, which should work particularly well in these circumstances, as the law of large numbers applies. This effect is slightly stronger for smaller sample sizes, when $\proc_n$ is quite not so close to its limit $\theta$ yet, which is also in agreement with our interpretation of $\optimal{\dec}$ being a compromise between the plug-in and \ac{DRO} method. In analogy, a more dominant penalization via an increased $\rate$ is favorable for smaller sample size $n$ and becomes a nuisance as $\proc_n$ approaches $\theta$.

\subsubsection{Asymptotic robust entropic risk measure}\label{sssec:rob:asym:ent:risk}

Let us revisit \Cref{ssec:LD:uncertainty}, i.e., we consider the sublinear expectation $\subexpecsub{\theta}{\cdot} =\sup_{\utheta\in B_\theta(R)}\expecsub{\utheta}{\cdot}$ with the ball $B_\theta(R)\subseteq\Theta\subseteq\Delta_d$.
Further, recall $\seq{\xi}{k}$ to be a stochastic process with the finite state space $\Xi=\{1,\ldots,d\}$, which is i.i.d.~with respect to the set of distributions $B_\theta(R)$ parametrized by some unknown $\theta\in\Theta$, i.e., for any $\utheta\in B_\theta(R)$ the process $\seq{\xi}{k}$ is i.i.d.~with $\probsub{\utheta}{\xi_k = i}= \utheta_i$ for $i=1,\ldots,d$. As outlined in \Cref{ssec:LD:uncertainty}, if $\seq{\proc}{n}$ is still the empirical measure process with $\sspa=\Delta_d$, and if $\theta\mapsto B_\theta(R)$ is a continuous\footnote{See \Cref{def_semicts_set_map} for the characterization of a set-valued continuous map.} and compact valued, then the asymptotic robust entropic risk measure $\phi^{\mathsf{RENT}}_\theta$ with an appropriately chosen $\Theta$ constitutes the \ac{LP}~\eqref{eq_LP2} with the robust relative entropy rate function $\rent$ given by \eqref{eq_def_robust_relative_entropy}.
Also here, this rate function can be shown to satisfy \Cref{ass_ratefct}, so that \Cref{thm_feas_dom_main} (and \Cref{thm_feas_dom_main_g}) apply.
\begin{proposition}\label{prop_robust_relentropy_satisfies}
    The robust relative entropy $ \rent$ defined in \eqref{eq_def_robust_relative_entropy} with $\theta\mapsto B_\theta(R)$ being a continuous set-valued map such that $B_\theta(R)$ is compact and non-empty satisfies \Cref{ass_ratefct}
    with $\sspa=\Delta_d$ and some compact $\Theta\subseteq\Delta_d\cap(0,1]^d$.
\end{proposition}
The proof of \Cref{prop_robust_relentropy_satisfies} is straightforward, as $\rent$ is continuous by \Cref{prop_relentropy_satisfies} and Berge's maximum theorem (\Cref{thm_BergeMaximum}), so that \Cref{ass_ratefct} follows immediately.\\

This new modeling framework still accommodates the sampling scheme of the classical newsvendor problem. However, the associated decision rule, based on a robust relative entropy rate function, is likely to be overly conservative, as it incorporates an additional layer of uncertainty via the sublinear expectation $\sup_{\utheta \in B_\theta(R)} \expecsub{\utheta}{\cdot}$.
Importantly, the framework can also capture \emph{aleatoric} uncertainty, such as data that is inherently noisy. This type of uncertainty is fundamentally different from the \emph{epistemic} uncertainty present in the classical newsvendor setting, where the uncertainty arises from incomplete knowledge of the true data-generating distribution $\theta$ and can be reduced by collecting more data.
In what follows, we consider a modified newsvendor problem under distributional (aleatoric) uncertainty, thereby illustrating the broader generality and modeling capabilities of our generalized framework introduced in Section~\ref{sec_other_opt_criteria}.

\paragraph{Robust newsvendor problem.} 
This time our newsvendor is again presented with the task of placing an order of a continuous amount of an underlying good under uncertain demand, which leads to a minimal expected loss; however, during the decision making process she only has data generated by multiple and slightly different sources at her disposal. Such a situation might occur if a new shop is to be opened, so that only observations from businesses with similar characteristics and circumstances can be used, as historical data is not yet available \cite[Section~6.2]{ref:rychener-2024}. In more precise terms, the newsvendor still has to minimize the regret of his decision with respect to an unknown and to her shop corresponding $\theta$ while dealing with sample paths, of which each was independently generated according to some $\utheta$ from a ball around $\theta$.

In a synthetic example, we select two distribution vectors $\utheta_1,\utheta_2\in B_{\theta}(R) = \{ \theta' \in \Theta \mid \norm{\theta' - \theta}_2 \leq R \}$, where $B_{\theta}(R)$ is the $L_2$ closed ball around $\theta$ with radius $R$ intersected with $\Theta$. Then, we generate $K$ i.i.d.~samples of length $N$ from the distribution $\utheta_1$ or $\utheta_2$ on $\Xi$, so that $\utheta_1$ is chosen for the $k$-th sample with probability $p_{\utheta_1}\in(0,1)$ and $\utheta_2$ otherwise.
This yields $K$ sample paths of the empirical distribution process $\seq{\proc}{n}$ of length $N$, denoted as $(z_1^{(k)},\ldots,z_N^{(k)})$ with $k=1,\ldots,K$. In the second row of \Cref{fig_plot_newsvendor_average_regrets}, the average regret of the decision $\optimal{\dec}$ 
under the unknown true distribution $\theta$ based on a sample of size $n=10,20,30,,\dots,N$ across $K$ experiments is plotted for $g\in\{c,\regret\}$ and for different values of $\beta$ and $\rate$.

By construction of our experiment, $\seq{\proc}{n}$ clearly cannot converge to $\theta$ and the (non-robust) law of large numbers does not hold, since any simulated path converges either to $\utheta_1\neq\theta$ or $\utheta_2\neq\theta$. Thus, the average regret of plug-in decision saturates at a value well above 0.
In that spirit, the distributional uncertainty has to be reflected in an adjustment of the notion of consistency, i.e., we consider a consistency set rather than consistency with respect to a single element of $\Theta$.

Our method shows the ability to adapt to such consequences of distributional (aleatoric) uncertainty in our data with both $g=c$ and $g=\regret$, where the latter gives much more significant signs of out-performance. 
However, the influence of parameters $\beta$ and $\rate$ is not as straightforward as above anymore. An increase in $\beta$ can help to decrease regret of the resulting decision even further and does not necessarily lead to an approximation of the plug-in decision. Larger $\rate$ on the other hand seems to destabilize the performance.

Although not displayed in \Cref{fig_plot_newsvendor_average_regrets}, the decision optimal under the standard Newsvendor setting, i.e.,  $\optimal{\dec}$ with $\ratefct=\ent$ and $g=\regret$, does also outperform $\decplug$ for $g=\regret$ and $\beta=0.5$; hence, this can be achieved provided that $\beta$ is sufficiently large for $\Gamma(\theta)\approx0$, but also small enough to provide flexibility, which differentiates it from the plug-in decision. This is again in line with our interpretation of $\optimal{\dec}$ compromising between  the plug-in and \ac{DRO} method (\Cref{subsec_problem_statement_consistency}).

\subsection{Empirical expectation of i.i.d.~process}\label{subsec_example_empmean_iid}
Cramér's theorem~\cite[Theorem 2.2.30]{iB1} allows to systematically identify large deviation rate functions corresponding to empirical means for i.i.d.~random variables. Whether the state space is some finite alphabet or a finite-dimensional real vector space, Cramér's theorem shows
the rate function to be the \defformat{Fenchel-Legendre transform} of the \defformat{logarithmic moment generating function} associated with the distribution of the underlying i.i.d.~process, which becomes practically relevant, once it can be computed numerically as a function in both $z$ and $\theta$.
For some distributions the Fenchel-Legendre transform is even available in closed form, which is the case in the following Gaussian example.

Let $\seq{\xi}{k}$ be an i.i.d.~stochastic process with $\xi_{k}\sim\mathcal{N}(\theta,\Sigma)$, where $\Sigma\in\real^{d\times d}$ is a known invertible covariance matrix and $\theta\in\Theta$ is an unknown mean vector from some compact $\Theta\subseteq\real^d$.
Then, according to Cram\'er's theorem the family of distributions $\seq{\law}{n}$ associated with the empirical mean process $\seq{\proc}{n}$, $\proc_n = \frac{1}{n}\sum_{k=1}^{n}\xi_{k}$ with state space $\sspa=\real^d$, satisfies an \ac{LDP}  with the rate function 
$$
\ratefct_{\theta}(z) = \frac{1}{2}(\theta-z)^{\top}\Sigma^{-1}(\theta-z),\quad (\theta,z)\in\Theta\times\sspa,
$$
as justified in \citep[Section~5]{iP1}. Hence, we have the \ac{LP}~\eqref{eq_LP2} for the asymptotic entropic risk measure $\phi^{\mathsf{ENT}}_{\theta}$ with the rate function $\ratefct_\theta$, which as a continuous function clearly satisfies \Cref{ass_ratefct}. Moreover, by the strong law of large numbers it holds $\proc_{n} \to \theta=\proc_\infty(\theta)$ $\bbP_\theta$-a.s.~as $n\to\infty$ and \Cref{ass_convergent_process} is satisfied.
Again this is an instance of the setting of \Cref{subsubsec:LD} but the decision is different from the finite-state i.i.d.~setting considered in \Cref{ex_plugin_via_LLN}.

It is straightforward to derive a number of other than Gaussian distributions with their respective rate functions in closed form, see \citep[Table~1]{iP1}. 
Notably, the Gärtner-Ellis theorem \cite[Section 2.3]{iB1} provides an extension to non-i.i.d.~cases and provides grounds for the Markov chain example below.

\paragraph{Portfolio optimization problem.} 
We consider a Markowitz-type portfolio selection problem with i.i.d.~Gaussian asset returns. More precisely, we consider asset allocation with or without short-selling represented by the decision vector of asset weights from $\dspa = \{x\in[-1,1]^d \mid \sum_{i=1}^d x_i = 1\}$ or $\dspa = \Delta_d$ respectively. The optimal portfolio shall minimize the negative of the standard mean-variance objective
\begin{equation}\label{eq:meam:variance:portfolio}
    c(x,\theta) = \mathbb{E}_\theta[-x^\top \xi] + \rho 
    \mathbb{V}_\theta[x^\top \xi] = -x^\top \theta + \rho x^\top \Sigma x,\quad (x,\theta)\in\dspa\times\Theta,
\end{equation}
where $\xi\in\real^d$ is the vector of asset returns that is described by a Gaussian distribution with unknown mean $\theta\in\Theta=\{\theta \in \real^d\mid\norm{\theta}\leq \delta\}$ for some $\delta>0$.

\begin{figure}[H]
    \begin{center}
        \includegraphics[width=0.33\textwidth]{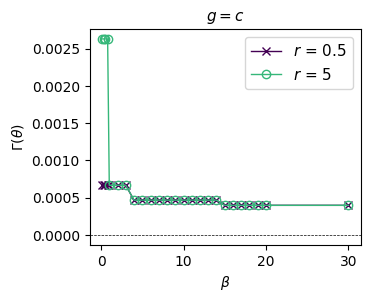}
    \end{center}
    \caption{Plot of the consistency gap of $\optimal{\dec}$ as defined by \eqref{eq_solution_star_g} with $\ratefct_\theta(z)=\frac{1}{2}(\theta-z)^{\top}\Sigma^{-1}(\theta-z)$, $(\theta,z)\in\Theta\times\sspa$, $g=c$ and $\rate\in\{0.5,5\}$ viewed as a function of the penalty function parameter $\beta$ within the portfolio optimization problem. We display consistency gap at the true $\theta$, i.e., $\Gamma(\theta)=\regret(\optimal{\dec}(\proc_\infty(\theta)),\theta)$. The cost \eqref{eq:meam:variance:portfolio} with $\rho = 1$ and data generating process parameters $d = 3$, $\Theta = \{\theta \in \real^3\mid\norm{\theta}_2\leq 10000 \}$, $\theta = (-0.2, 0.6,0.35)^\top$ and covariance matrix $\Sigma = [(2.819, 1.726, 1.917), (1.726, 1.297, 1.081), (1.917, 1.081, 2.717)]$ are considered.
    }
    \label{fig_plot_Gaussian_GapBeta}
\end{figure}

Overall, the simulation results for the portfolio optimization problem are not much different to those provided for the newsvendor model above. This is particularly true about the analysis of the consistency gap at an unknown $\theta$, $\Gamma(\theta)$, which converges to zero as $\beta\to\infty$. To complete the heuristics discussed in \Cref{subsec_problem_statement_consistency}, where this Gaussian setup proved to be useful, we provide the corresponding empirical findings in \Cref{fig_plot_Gaussian_GapBeta}.

\subsection{Finite state Markov chain - pair empirical measure}\label{subsec_example_finitestate_MC_pair}

Turning our attention to setups not requiring the (robust) i.i.d. assumption, we will adapt the example of \cite{iP1} and \cite[Section 3.1.3]{iB1} to our setting. Let $\seq{\xi}{k}$ be a time-homogeneous ergodic Markov chain taking values in the finite state space $\Xi\defeq\{1,\ldots,d\}$
with some deterministic initial state $\xi_0\in\Xi$ and $\lim_{k\to\infty}\probsub{\theta}{\xi_{k} = i, \xi_{k+1}=j} =\theta_{ij}>0$ for all $i,j\in\Xi$, so that the matrix $\theta\in\Theta$ encodes the stationary distribution of the pair $(\xi_k,\xi_{k+1})$. Hence, $\theta$ must be implicitly such that $\theta\in(0,1]^{d\times d}$, $\sum_{i,j\in\Xi}\theta_{ij} = 1$ and $\sum_{j\in\Xi}\theta_{ij}=\sum_{j\in\Xi}\theta_{ji}$ for all $i\in\Xi$. In alignment with our overall assumptions, $\Theta$ might then be a set of all $\theta\in[\epsilon,1]^{d\times d}$ with such properties for some $0<\epsilon\ll1$.
Then, the family of distributions $\seq{\law}{n}$ associated with the process of pair empirical measures $\seq{\proc}{n}$, $(\proc_n)_{ij} = \frac{1}{n}\sum_{k=1}^{n}\indfctset{(\xi_{k-1},\xi_k)=(i,j)}$ for $i,j\in\Xi$, with state space $\sspa = \Delta_{d^2}$ satisfies an \ac{LDP} with the rate function
being the so called \defformat{conditional relative entropy}
\begin{equation}\label{eq_def_cond_relentropy}
    \ratefct_{\theta}(z) =\sum_{i,j\in\Xi} z_{ij} \left( \log \frac{z_{ij}}{\sum_{k\in\Xi} z_{ik}} - \log \frac{\theta_{ij}}{\sum_{k\in\Xi} \theta_{ik}}\right),\quad(\theta,z)\in\Theta\times\sspa,
\end{equation}
with the conventions $0\log0=0$ and $0\log\frac{0}{0}=0$. This rate function is continuous on $\Theta\times\sspa$ (\cite{iP1}) and the \Cref{ass_ratefct} is satisfied.

The \Cref{ass_convergent_process} holds by the ergodic theorem for Markov chains, which ensures that $\proc_{n} \to \theta=\proc_\infty(\theta)$ $\bbP_\theta$-a.s. as $n\to\infty$. With respect to consistency however, the rate function \eqref{eq_def_cond_relentropy} lacks one of the desired properties listed in \Cref{subsec_problem_statement_consistency}, which in fact does hold for the relative entropy. Enhancing $\optimal{\dec}$ with the consistency properties of the plug-in estimator \eqref{eq_plug} by choosing a large penalty function parameter $\beta$ is not easily achievable here, as $\ratefct_{\utheta}(\theta)=0$ for some $\utheta,\theta\in\Theta$ does not necessarily imply $\utheta=\theta$.
The resulting instability in the relationship between the consistency gap $\Gamma$ and $\beta$ is well visible in our simulations.
\section{Proofs}\label{sec_proofs}

To derive the proofs of \Cref{thm_feas_dom_main} and other results, we only need to treat the generalized setting \eqref{eq_optp_main_g} with the function $g\pcolon \dspa\times\Theta\map\real$ defined in \Cref{sec_other_opt_criteria}. The modified problem \eqref{eq_optp_modified_infinf_g} requires a separate proof, as an application of the Sion's minimax theorem~\cite{Sion:1958} is necessary.
The functions $\optimal{\dec}$, $\optimal{\off}$ and $\supfct$ are defined in \eqref{eq_solution_star_g} and \eqref{eq_G_g}.

\begin{lemma}[Continuity of $\supfct$]\label{lemma_optimal_G_g_cts}
    Let $\ratefct\pcolon \Theta\times\sspa\map[0,\infty]$ satisfy \Cref{ass_ratefct}. Then, the function $\supfct\pcolon\dspa\times\sspa\to\real$ defined by \eqref{eq_G_g}
    is continuous on $\dspa\times\sspa$ and the maximum is attained. Moreover, $x\mapsto\supfct (x,z)$ is strictly convex for all $z\in\mathcal{Z}$.
\end{lemma}
\begin{proof}
To show that $-\infty < \supfct < \infty$, let $x \in \dspa$, $z \in \sspa$ be fixed.  
By \Cref{ass_ratefct}~\ref{item_ass_ratefct_dom}, there exists $\theta_z \in \Theta$ such that $0 \leq \ratefct_{\theta_z}(z) < \infty$. It follows that
\[
\supfct(x, z) \geq g(x, \theta_z) - \inv{\ell}(\ratefct_{\theta_z}(z) - \rate) > -\infty.
\]
Moreover, since $g$ is continuous and $\Theta$ is compact, we have $\sup_{\theta \in \Theta} g(x, \theta) < \infty$. Also, $\inv{\ell}$ is increasing and $\ratefct_{\theta}(z) - \rate \geq -\rate$ for all $\theta \in \Theta$, so that
\[
\supfct(x, z) \leq \max_{\theta \in \Theta} g(x, \theta) - \inv{\ell}(-\rate) < \infty.
\]

We now show that $\supfct$ is continuous.
First, that $\supfct$ is \ac{usc} follows directly from \Cref{ass_ratefct}~\ref{item_ass_ratefct_lsc}. Indeed, the joint lower semi-continuity of the rate function implies via \Cref{lemma_penalty_lsc_composit_lsc} that the function
\[
(x, z, \theta) \mapsto g(x, \theta) - \inv{\ell}(\ratefct_{\theta}(z) - \rate)
\]
is \ac{usc}. It follows that the supremum in the definition of $\supfct$ is attained over the compact set $\Theta$, i.e., the supremum is in fact a maximum. Moreover,
 for any sequence $(x_n, z_n)_{n \in \mathbb{N}}$ converging to $(x, z)$, there exist a sequence of maximizers $(\theta_n)_{n \in \mathbb{N}}$ for $G(x_n,z_n)$ and some $\theta\in\Theta$ with $\theta_n \to \theta$ (along a subsequence), such that
\[
\begin{aligned}
\supfct(x, z) &\geq g(x, \theta) - \inv{\ell}(\ratefct_{\theta}(z) - \rate) \\
&\geq \limsup_{n \to \infty} \left[ g(x_n, \theta_n) - \inv{\ell}(\ratefct_{\theta_n}(z_n) - \rate) \right] = \limsup_{n \to \infty} \supfct(x_n, z_n).
\end{aligned}
\]
Second, we show that $\supfct$ is also \ac{lsc}. Let $(x_n, z_n)_{n \in \mathbb{N}}$ converging to $(x, z)$ and $\theta\in\Theta$ be a maximizer for $G(x,z)$ so that $\ratefct_\theta(z) < \infty$. By \Cref{ass_ratefct}~\ref{item_ass_ratefct_edgects}, there exists a sequence $(\theta_n)_{n \in \mathbb{N}}$ with $\theta_n \to \theta$ and $\ratefct_{\theta_n}(z_n) \to \ratefct_\theta(z)$ such that
\[
\begin{aligned}
\supfct(x, z) &= g(x, \theta) - \inv{\ell}(\ratefct_\theta(z) - \rate) \\
&= \lim_{n \to \infty} \left[ g(x_n, \theta_n) - \inv{\ell}(\ratefct_{\theta_n}(z_n) - \rate) \right] \\
&\leq \liminf_{n \to \infty} \maxb{\theta \in \Theta}{g(x_n, \theta) - \inv{\ell}(\ratefct_{\theta}(z_n) - \rate)} = \liminf_{n \to \infty} \supfct(x_n, z_n).
\end{aligned}
\]

Finally, the strict convexity of $x \mapsto \supfct(x, z)$ follows since $\supfct(x, z)$ is the pointwise maximum over $\theta \in \set{\theta\in\Theta}{\ratefct_\theta(z)<\infty}$ of the strictly convex functions
$x \mapsto g(x, \theta) - \inv{\ell}(\ratefct_\theta(z) - \rate)$.
As the pointwise maximum of strictly convex functions is again strictly convex, the claim follows.
\end{proof}

\begin{lemma}[Continuity of $\optimal{\dec}$]\label{lemma_optimal_X_g_cts}
    Under the \Cref{ass_ratefct} the decision function $\optimal{\dec}\pcolon\sspa\to\dspa$ defined by \eqref{eq_solution_star_g} is continuous on $\sspa$.
\end{lemma}

\begin{proof}
By \Cref{lemma_optimal_G_g_cts}, the function $\supfct$ is jointly continuous. Hence, we can apply Berge's Maximum~\Cref{thm_BergeMaximum} (with the constant set-valued map $\sspa \ni z \mapsto \dspa$) to conclude that the set-valued map
\begin{equation}\label{eq_argmin_map_g}
    \sspa \ni z \mapsto \set{x^\star \in \dspa}{x^\star \in \Argmin_{x \in \dspa} \supfct(x, z)} \in \mathcal{P}(\dspa)
\end{equation}
is \ac{usc} (in the sense of \Cref{def_semicts_set_map}) and has non-empty, compact values. By \Cref{lemma_optimal_G_g_cts}, the function $x \mapsto \supfct(x, z)$ is strictly convex for each $z \in \sspa$, so the mapping in~\eqref{eq_argmin_map_g} is single-valued and coincides by definition with $z \mapsto \optimal{\dec}(z)$.

Finally, since single-valued, \ac{usc} set-valued maps are continuous by \Cref{def_semicts_set_map}, it follows that $z \mapsto \optimal{\dec}(z)$ is continuous.
\end{proof}

\begin{proof}[\proofheading{Proof of \Cref{thm_feas_dom_main_g}}] 
By \Cref{lemma_optimal_G_g_cts} and \Cref{lemma_optimal_X_g_cts}, the function $G$ is jointly continuous, the argmin in~\eqref{eq_optp_main_g} is single-valued, and $(\optimal{\dec}, \optimal{\off}) \in \decset \times \uset$.

\begin{enumerate}[label=\enumthm]
    \item \label{item_proof_optim_main_a_g}
    Fix $\theta \in \Theta$. Since the function $\sspa \to \real$ given by $z\mapsto \ell(g(\optimal{\dec}(z), \theta) - \optimal{\off}(z))$
    is continuous, it follows from the Laplace principle~\eqref{eq_LP2} that
    \begin{align}
        &\phi_\theta\left[\left(\ell\left(g(\optimal{\dec}(\proc_n), \theta) - \optimal{\off}(\proc_n)\right)\right)_{n \in \nat}\right]  \nonumber \\
        &\qquad = \sup_{z \in \sspa} \left\{ \ell\left(g(\optimal{\dec}(z), \theta) - \supfct(\optimal{\dec}(z), z)\right) - \ratefct_\theta(z) \right\} \phantom{\Big\{ } \nonumber \\
        &\qquad \leq \sup_{\substack{z \in \sspa,\, \ratefct_\theta(z) < \infty}} \left\{ \ell\Big( \maxb{x \in \dspa}{g(x, \theta) - \supfct(x, z)} \Big) - \ratefct_\theta(z) \right\}. \label{eq_proof_feasibility_main_onethmarg_g} 
    \end{align}
    Moreover, for every $z \in \sspa$, it holds that
    \begin{align*}
        \maxb{x \in \dspa}{g(x, \theta) - \supfct(x, z)}
        &= \maxb{x \in \dspa}{g(x, \theta) - \maxb{\utheta \in \Theta}{g(x, \utheta) - \inv{\ell}(\ratefct_\utheta(z) - \rate)}} \\
        &\leq \maxb{x \in \dspa}{g(x, \theta) - g(x, \theta) + \inv{\ell}(\ratefct_\theta(z) - \rate)} \\
        &= \inv{\ell}(\ratefct_\theta(z) - \rate),
    \end{align*}
    which, together with~\eqref{eq_proof_feasibility_main_onethmarg_g} and using the fact that $\inv{\ell}$ is the inverse of $\ell$, implies
    \begin{align*}
        \phi_\theta\left[\left(\ell\left(g(\optimal{\dec}(\proc_n), \theta) - \optimal{\off}(\proc_n)\right)\right)_{n \in \nat}\right]
        &\leq \sup_{\substack{z \in \sspa,\, \ratefct_\theta(z) < \infty}} \left\{ \ell\left(\inv{\ell}(\ratefct_\theta(z) - \rate)\right) - \ratefct_\theta(z) \right\} \leq -\rate,
    \end{align*}
    which completes the proof of assertion~\ref{item_proof_optim_main_a_g}.

    \item Let $(\dec, \off)$ be a feasible pair in~\eqref{eq_optp_main_g}. Since $z \mapsto \ell\left(g(\dec(z), \theta) - \off(z)\right)$ is continuous for all $\theta \in \Theta$, it follows from the Laplace principle~\eqref{eq_LP2} that
    \[
    \sup_{z \in \sspa} \left\{ \ell\left(g(\dec(z), \theta) - \off(z)\right) - \ratefct_\theta(z) \right\} = \phi_{\theta}\left[\ell\left(g(\dec(\proc_n), \theta) - \off(\proc_n)\right)_{n \in \nat}\right] \leq -\rate,
    \]
    so that $\ell\left(g(\dec(z), \theta) - \off(z)\right) \leq \ratefct_\theta(z) - \rate$ for all $\theta \in \Theta$ and $z \in \sspa$. Since $\inv{\ell}$ is increasing, this further implies
    \[
    \inv{\ell}(\ratefct_\theta(z) - \rate) \geq \inv{\ell}\left(\ell\left(g(\dec(z), \theta) - \off(z)\right)\right) = g(\dec(z), \theta) - \off(z) \quad \forall \theta \in \Theta,\, \forall z \in \sspa.
    \]
    Rearranging the previous inequality and taking the supremum over $\theta \in \Theta$ yields
    \begin{align*}
        \off(z) &\geq \supb{\theta \in \Theta}{g(\dec(z), \theta) - \inv{\ell}(\ratefct_\theta(z) - \rate)} \\
        &\geq \min_{x \in \dspa} \supb{\theta \in \Theta}{g(x, \theta) - \inv{\ell}(\ratefct_\theta(z) - \rate)} \\
        &= \min_{x \in \dspa} \supfct(x, z) = \optimal{\off}(z)
    \end{align*}
    for all $z \in \sspa$. This shows that $(\optimal{\dec},\optimal{\off})$ is optimal.\qedhere
\end{enumerate}
\end{proof}

\begin{proof}[\proofheading{Proof of \Cref{prop_feas_dom_modified_infinf_g}}] 
We show the two assertions separately. 
\begin{enumerate}[label=\enumproposition]
    \item Once again, $\supfct \in \costset$ follows immediately from \Cref{lemma_optimal_G_g_cts}. Since the mapping $z \mapsto \min_{x \in \dspa} g(x, \theta) - \min_{x \in \dspa} \supfct(x, z)$ is continuous by Berge's maximum theorem (\Cref{thm_BergeMaximum}), it follows from the Laplace principle~\eqref{eq_LP2} that
    \begin{align*}
        &\phi_{\theta}\left[\left(\ell\Big( \min_{x \in \dspa} g(x, \theta) - \min_{x \in \dspa} \supfct(x, \proc_n) \Big)\right)_{n \in \nat} \right] \\
        &\qquad = \sup_{z \in \sspa} \left\{ \ell\Big( \min_{x \in \dspa} g(x, \theta) - \min_{x \in \dspa} \supfct(x, z) \Big) - \ratefct_\theta(z) \right\} \\
        &\qquad \leq \sup_{\substack{z \in \sspa,\, \ratefct_\theta(z) < \infty}} \left\{ \ell\Big( \maxb{x \in \dspa}{g(x, \theta) - \supfct(x, z)} \Big) - \ratefct_\theta(z) \right\},
    \end{align*}
    which holds for all $\theta \in \Theta$, where the inequality follows from \Cref{lemma_inf_minus_inf_leq_sup}. The remainder of the argument is analogous to the proof of \Cref{thm_feas_dom_main_g}~\ref{item_thm_feas_dom_main_feas_g}.

    \item Let $\modsupfct \in \costset$ be a feasible solution to~\eqref{eq_optp_modified_infinf_g}. Again by Berge's maximum theorem (\Cref{thm_BergeMaximum}), the mapping $z \mapsto \ell\left(\min_{x \in \dspa} g(x, \theta) - \min_{x \in \dspa} \modsupfct(x, z)\right)$ is continuous. Hence, for every $\theta \in \Theta$, the Laplace principle~\eqref{eq_LP2} yields
    \begin{align*}
        &\sup_{z \in \sspa} \left\{ \ell\left(\min_{x \in \dspa} g(x, \theta) - \min_{x \in \dspa} \modsupfct(x, z)\right) - \ratefct_\theta(z) \right\} \\
        &\qquad= \phi_{\theta}\left[\left( \ell\Big( \min_{x \in \dspa} g(x, \theta) - \min_{x \in \dspa} \modsupfct(x, \proc_n) \Big) \right)_{n \in \nat} \right] \leq -\rate,
    \end{align*}
    which implies that
    \[
    \ell\left( \min_{x \in \dspa} g(x, \theta) - \min_{x \in \dspa} \modsupfct(x, z) \right) \leq \ratefct_\theta(z) - \rate
    \quad \forall \theta \in \Theta,\ \forall z \in \sspa.
    \]
    Since $\inv{\ell}$ is increasing, it follows that
    \[
    \inv{\ell}(\ratefct_\theta(z) - \rate)
    \geq \inv{\ell}\left( \ell\left( \min_{x \in \dspa} g(x, \theta) - \min_{x \in \dspa} \modsupfct(x, z) \right) \right)
    = \min_{x \in \dspa} g(x, \theta) - \min_{x \in \dspa} \modsupfct(x, z)
    \]
    for all $\theta \in \Theta$ and $z \in \sspa$. Rearranging yields
    \[
    \min_{x \in \dspa} \modsupfct(x, z)
    \geq \supb{\theta \in \Theta}{\min_{x \in \dspa} g(x, \theta) - \inv{\ell}(\ratefct_\theta(z) - \rate)},
    \quad \forall z \in \sspa.
    \]
    
    Using the concavity assumptions, it follows from Sion's minimax theorem~\cite{Sion:1958}, as discussed in \Cref{rem_concave}, that for every $z \in \sspa$,
    \[
    \sup_{\theta \in \Theta} \min_{x \in \dspa} \left\{ g(x, \theta) - \inv{\ell}(\ratefct_\theta(z) - \rate) \right\}
    = \min_{x \in \dspa} \max_{\theta \in \Theta} \left\{ g(x, \theta) - \inv{\ell}(\ratefct_\theta(z) - \rate) \right\}
    = \min_{x \in \dspa} \supfct(x, z),
    \]
    which, in combination with the previous inequality, implies that
    \[
    \min_{x \in \dspa} \supfct(x, z) \leq \min_{x \in \dspa} \modsupfct(x, z),
    \quad \forall z \in \sspa. \qedhere
    \]
\end{enumerate}
\end{proof}

\paragraph{Acknowledgments.} The authors would like to thank Daniel Kuhn, Ariel Neufeld and Jan Beran for helpful discussions on an early version of this paper.
\paragraph{Supplements.} Code for the experiments presented in  \Cref{sec_examples} is available at \url{https://github.com/radeksalac/Asymptotic-Optimality-in-Data-Driven-Decision-Making}

\appendix

\section*{Appendix}\label{sec_appendix} 
\renewcommand{\theequation}{A.\arabic{equation}}  
\renewcommand{\thefigure}{A.\arabic{figure}}      
\renewcommand{\thetable}{A.\arabic{table}}        

\renewcommand{\thethm}{A.\arabic{thm}}

Let $X$ and $Y$ be metric spaces. We denote by $\mathcal{P}(Y)$ the power set of $Y$, i.e., the collection of all subsets of $Y$.
Further, any function $f\pcolon X\map\realinf$ is \acf{lsc} if and only if for every $x_0\in X$,
$$\liminf_{x\converge x_0} f(x)\geq f(x_0),$$
or equivalently, if for each $y\in\real$, the sub-level set $\set{x\in X}{f(x)\leq y}$ is closed in $X$.
It is \acf{usc} if and only if for every $x_0\in X$,
$$\limsup_{x\converge x_0} f(x)\leq f(x_0),$$
or equivalently, if for each $y\in\real$, the super-level set $\set{x\in X}{f(x)\geq y}$ is closed in $X$.

We recall that a penalty function $\ell\pcolon \real \to \real$ can be uniquely extended to a continuous function $\ell\pcolon \realinf \to \realinf$ by setting $\ell(-\infty) \defeq -\infty$ and $\ell(\infty) \defeq \infty$. We frequently work with this extended version when necessary.

\begin{lemma}\label{lemma_penalty_lsc_composit_lsc}
Let $g\pcolon X\map\realinf $ be \ac{lsc} and $\ell\pcolon\realinf\map\realinf$ be a penalty function.
    Then, $\ell\circ g\pcolon X\map\realinf $ is also \ac{lsc}.
\end{lemma}
\begin{proof} For every $y\in \real$, it holds  $\set{x\in X}{\ell\circ g(x)\leq y}=\set{x\in X}{g(x)\leq \inv{\ell}(y)}$, which is a closed set, because $g$ is \ac{lsc} and $\inv{\ell}(y)\in\real$. 
\end{proof}

\begin{definition}[Semi-continuous set-valued maps]\label{def_semicts_set_map}
    Let $\setmap\pcolon  X\map\mathcal{P}(Y)$ be a set-valued map. 
    The mapping $\setmap$ is called \acf{lsc}, if for each $x_0\in X$ and every open set $G\subseteq Y$ with $\setmap (x_0) \cap G\neq\emptyset$, there exists a neighborhood $U\subseteq X$ of $x_0$ such that $\setmap (x) \cap G\neq\emptyset$ for all $x\in U$. Moreover, the mapping $\setmap$ is called \acf{usc}, if for each $x_0\in X$
    and every open set $G\subseteq Y$ with $\setmap (x_0) \subseteq G$, there exists a neighborhood $U\subseteq X$ of $x_0$ such that $\setmap (x) \subseteq G$ for all $x\in U$.
    Finally, $\setmap$ is called continuous, if it is both \ac{usc} and \ac{lsc}.
\end{definition}

As a number of our results rely on the Berge's maximum theorem~\cite{berge1997topological,iB17}, we display its standard version in the following theorem.
\begin{thm}[Berge's maximum theorem]\label{thm_BergeMaximum}
Let $\varphi\pcolon X\times Y\map \real$ be a jointly continuous function and $\setmap\pcolon  X\map\mathcal{P}(Y)$ be a continuous set-valued map such that $\setmap (x)\subseteq Y$ is compact and non-empty for all $x\in X$. Then, it holds:
    \begin{enumerate}[label=\enumthm]
        \item The function $M\pcolon X\map\real$ defined by 
        $M(x)\defeq\sup\set{\varphi(x,y)}{y\in\setmap (x)}$
        is continuous and the supremum is attained for all $x\in X$.
        \item The ``\textit{argmax}'' set-valued map $X\to\mathcal{P}(Y)$, 
        \begin{equation*}
            x \mapsto \set{y\in\setmap (x)}{\varphi(x,y)=M(x)}
        \end{equation*}
        has non-empty compact values and is \ac{usc}.
        \label{item_b_thm_BergeMaximum}
    \end{enumerate}
\end{thm}

Besides the statements about (semi-)continuity of set-valued maps and functions to the extended real line, we also make use of the following simple rule for bounding the difference of two infima.

\begin{lemma}\label{lemma_inf_minus_inf_leq_sup}
    Let $f,g\pcolon X\map\real$ be so that either 
     $\inf_{x\in X} f(x)$ or $\inf_{x \in X} g(x)$ is real-valued.
    Then, it holds
    \begin{equation*}
        \inf_{x\in X} f(x) - \inf_{x \in X} g(x) \leq \sup_{x\in X} f(x) - g(x).
    \end{equation*}
\end{lemma}
\begin{proof}
    For every $\varepsilon>0$, there exists an $x_\varepsilon\in X$ such that $\inf_{x\in X} g(x) \vee \left(-\frac{1}{\varepsilon}\right) \geq g(x_\varepsilon) -\varepsilon $. Hence, 
        \[ \inf_{x\in X} f(x) - \left(\inf_{x \in X} g(x)\vee \left(-\frac{1}{\varepsilon}\right)\right)\leq f(x_\varepsilon)-g(x_\varepsilon) +\varepsilon\leq\sup_{x\in X} f(x) - g(x)+\varepsilon.\]
    The claim follows by letting $\varepsilon \to 0$.
\end{proof}

In the following, we provide proofs for supportive statements used in \Cref{subsec_example_empmeasure_iid}.
These include
\Cref{prop_relentropy_satisfies} and \Cref{prop_robust_relentropy_satisfies}, showing that \Cref{ass_ratefct} holds for the finite state i.i.d. settings described in \Cref{subsec_example_empmeasure_iid}.

Preparatory to the proof of \Cref{prop_relentropy_satisfies}, we first show the following lemma.
\begin{lemma}\label{lemma_relentropy_discr_jointlyctseffdom}
The relative entropy $\ent$ defined in \eqref{eq_def_relative_entropy} constrained to its effective domain is jointly continuous, i.e., for all $(\theta,z),(\theta_n,z_n)\in\dom_{\Delta_d\times\Delta_d} \ent =\set{(\theta',z')\in\Delta_d\times\Delta_d}{z' \ll \theta'}$ with $z_n\converge z$ and $\theta_n\converge\theta$, it holds $\ent_{\theta_n}(z_n)\converge\ent_\theta(z)$ as $n\converge\infty$.
\end{lemma}
\begin{proof}
    Let $(\theta,z),(\theta_n,z_n)\in\dom_{\Delta_d\times\Delta_d} \ent$ such that $z_n\converge z$ and $\theta_n\converge\theta$. It is enough to show that $$\lim_{n\converge\infty} \ (z_n)_i \log\frac{(z_n)_i}{(\theta_n)_i}=(z)_i \log\frac{(z)_i}{(\theta)_i}$$ for all $i\in\{1,\dots,d\}$.
    To do so, let $i\in\{1,\dots,d\}$ and consider the case $(\theta)_i,(z)_i>0$. Then, w.l.o.g., $(\theta_n)_i,(z_n)_i>0$ and the result follows by the fact that $(x,y)\mapsto x\log x/y$  is  continuous on $(0,1]^2$. Next, assume $(\theta)_i>0$ and $(z)_i=0$. Then, we have
    \begin{equation*}
        \abs{(z_n)_i\log\frac{(z_n)_i}{(\theta_n)_i} - (z)_i\log\frac{(z)_i}{(\theta)_i}} =
        \abs{(z_n)_i\log\frac{(z_n)_i}{(\theta_n)_i}}\leq
        \abs{(z_n)_i\log(z_n)_i} + \abs{(z_n)_i\log(\theta_n)_i},
    \end{equation*}
    where the last term converges to $0+0\abs{\log (\theta)_i} =0$ as $n\converge\infty$. Finally, we inspect the case $(\theta)_i=(z)_i=0$. It holds
    \begin{equation*}
        \abs{(z_n)_i\log\frac{(z_n)_i}{(\theta_n)_i} - (z)_i\log\frac{(z)_i}{(\theta)_i}} =
        \abs{(z_n)_i\log\frac{(z_n)_i}{(\theta_n)_i}}\leq(z_n)_i\frac{1-\frac{(z_n)_i}{(\theta_n)_i}}{\frac{(z_n)_i}{(\theta_n)_i}}=(\theta_n)_i- (z_n)_i\converge 0
    \end{equation*}
    as $n\converge\infty$, while using the fact that $\abs{\log x}\leq\frac{1-x}{x}$ for all $x\in[0,1]$.
\end{proof}

\begin{proof}[\proofheading{Proof of \Cref{prop_relentropy_satisfies}}]
    In the following, we only have to deal with the case $\Theta=\Delta_d=\sspa$, as for any compact $\Theta\subseteq\Delta_d\cap(0,1]^d$ the rate function $\ent$ is  finite by definition and even continuous by \Cref{lemma_relentropy_discr_jointlyctseffdom}.
    Let $z\in\Delta_d$ and take any $\theta\in\Delta_d\cap(0,1]^d$. Then, $\ent_\theta(z)<\infty$, which already proves \Cref{ass_ratefct}~\ref{item_ass_ratefct_dom}. Next, we address \Cref{ass_ratefct}~\ref{item_ass_ratefct_lsc}, i.e., we show that $\ent$ is \ac{lsc} on $\Theta\times\sspa$.
    Let $z_n,z\in\sspa$ and $\theta_n,\theta\in\Theta$ such that $z_n\converge z$ and $\theta_n\converge \theta$ as $n\converge\infty$.
    We need to show $\liminf_{n\converge\infty}\ent_{\theta_n}(z_n)\geq \ent_\theta(z)$, which holds trivially if $\ent_{\theta_n}(z_n)=\infty$ for all $n\in\nat$. Hence, w.l.o.g., let $(\ent_{\theta_n}(z_n))_{n\in\nat}$ be a real-valued sequence. If $\ent_\theta(z)<\infty$, then \Cref{lemma_relentropy_discr_jointlyctseffdom} applies. In the other case, $\ent_\theta(z)=\infty$, there is again an $i\in\{1,\dots,d\}$ such that $(\theta)_{i}=0$ and $(z)_{i}>0$. But then, for an arbitrary small $(z)_{i}>\varepsilon_i>0$, there is a $n_\varepsilon\in\nat$,so that for all $n\geq n_{\varepsilon_i}$, it holds
    \begin{equation*}
         (z_n)_i \log\frac{(z_n)_i}{(\theta_n)_i}  \geq  ((z)_i-\varepsilon_i) \log\frac{(z)_i-\varepsilon_i}{\varepsilon_i}.
    \end{equation*}
    Hence, letting $\varepsilon_i$ converge to 0 gives $\liminf_{n\converge\infty}\ (z_n)_i \log\frac{(z_n)_i}{(\theta_n)_i}=\infty$. Consequently, we have $\liminf_{n\converge\infty}\ent_{\theta_n}(z_n)=\infty=\ent_\theta(z)$, because 
    \begin{equation*}
        \liminf_{n\converge\infty} \ \ent_{\theta_n}(z_n) = \sum_{i\not\in A_{\theta,z}} \underbrace{\liminf_{n\converge\infty} \ (z_n)_i \log\frac{(z_n)_i}{(\theta_n)_i}}_{\in\real} + \sum_{i\in A_{\theta,z}} \underbrace{\liminf_{n\converge\infty} \ (z_n)_i \log\frac{(z_n)_i}{(\theta_n)_i}}_{=\infty}=\infty,
    \end{equation*}
    where $A_{\theta,z}\defeq\set{i\in\{1,\dots,d\}}{(\theta)_{i}=0\text{ and }(z)_{i}>0}$, and where the first sum is finite by the proof of \Cref{lemma_relentropy_discr_jointlyctseffdom}.

    Finally, we will show \Cref{ass_ratefct}~\ref{item_ass_ratefct_edgects}. Let $(\theta,z)\in\dom_{\Delta_d\times\Delta_d} \ent$ and $\seq{z}{n}\subseteq\sspa$ with $z_n\converge z$ as $n\converge\infty$. We have to show that there is a sequence $\seq{\theta}{n}$ with $\theta_n\to\theta$ and $ \ent_{\theta_n}(z_n)\to\ent_\theta(z)$ as $n\converge\infty$.
    If $\seq{z}{n}\subseteq\dom_\sspa \ent_\theta$, then we are done by choosing the constant sequence $\seq{\theta}{n}$, $\theta_n\defeq\theta$, as $\ent_\theta$ is continuous and finite on $\dom_\sspa\ent_\theta$, see for example \cite[p.~13]{iB1}. In the other case, $\seq{z}{n}\not\subseteq\dom_\sspa \ent_\theta$, we assume w.l.o.g. that $\ent_\theta(z_n)=\infty$ for all $n\in\nat$. Then, for all $n$ there is $i_n\in\{1,\dots,d\}$ such that $(\theta)_{i_n}=0$ and $(z_n)_{i_n}>0$. We define the desired sequence via 
    $\theta_n \defeq \frac{n-1}{n} \theta + \frac{1}{n} z_n$, which converges to $\theta$. 
    Then, $\ent_{\theta^n}(z^n) < \infty$ for all $n \in \mathbb{N}$, and the claim follows from \Cref{lemma_relentropy_discr_jointlyctseffdom}.
\end{proof} 

Finally, we prove the finite state Sanov's theorem adapted to our robust setting from 
Section~\ref{ssec:LD:uncertainty}.
\begin{proposition}[Robust Sanov's theorem for finite alphabets]\label{prop_robust_finite_sanov}
Let $\seq{\xi}{k}$ be a stochastic process on the finite state space $\{1, \ldots, d\}$, which is i.i.d.\ with respect to a non-empty set of distributions $B_\theta(R) \subseteq \Theta \subseteq \Delta_d$, parametrized by $\theta \in \Theta$. That is, under any $\utheta \in B_{\theta}(R)$, the process $\seq{\xi}{k}$ is i.i.d.\ with $\probsub{\utheta}{\xi_k = i} = \utheta_i$ for $i = 1, \ldots, d$. Let $\rent$ be the robust relative entropy rate function given by~\eqref{eq_def_robust_relative_entropy}. Assume that $B_\theta(R)$ and $\Theta$ are such that $\rent_\theta \pcolon \sspa \to [0, \infty)$ is continuous. Further, $\proc_n\defeq \sum_{i=1}^d \left(\frac{1}{n}\sum_{k=1}^{n}\indfctset{\xi_k=i}\right)e_i$ defines the empirical measure process $\seq{\proc}{n}$ associated with $\seq{\xi}{k}$. 
    Then, for every set $A$ of probability vectors in $\Delta_d$, it holds
    \begin{align}\label{eq_robust_Sanov}
    \begin{split}
        -\inf_{z\in \interior{A}} \rent_\theta(z) &\leq \liminf_{n\to\infty}\frac{1}{n}\log\sup_{\utheta\in B_\theta(R)}\probsub{\utheta}{\proc_n\in A}\\
        &\leq\limsup_{n\to\infty}\frac{1}{n}\log\sup_{\utheta\in B_\theta(R)}\probsub{\utheta}{\proc_n\in A}\leq-\inf_{z\in A} \rent_\theta(z),
    \end{split}
    \end{align}
\end{proposition}

\begin{proof}
     We adapt the proof of \cite[Theorem~2.1.10]{iB1}. 
    Let $A\subseteq\Delta_d$, $\theta\in\Theta$ and define 
    $\mathcal{L}_n\defeq\set{z\in\sspa}{\exists y\in\{1,\ldots,d\}^n : z = \sum_{i=1}^d \left(\frac{1}{n}\sum_{k=1}^{n}\indfctset{y_k=i}\right)e_i}$.
    By taking the supremum over $B_\theta(R)$ on both sides of the inequalities (2.1.12) and (2.1.13) in the proof of \cite[Theorem~2.1.10]{iB1}, we get
    \begin{align*}
        \sup_{\utheta\in B_\theta(R)}\probsub{\utheta}{\proc_n\in A}&\leq
        (n+1)^{d}\sup_{\utheta\in B_\theta(R)}\exp\left(-n\inf_{z\in A\cap \mathcal{L}_n}\ent_{\utheta}(z)\right)\\
        &= (n+1)^{d}\exp\left(-n\inf_{z\in A\cap \mathcal{L}_n}\rent_{\theta}(z)\right)
    \end{align*}
    and
    \begin{align*}
        \sup_{\utheta\in B_\theta(R)}\probsub{\utheta}{\proc_n\in A}&\geq (n+1)^{-d}\sup_{\utheta\in B_\theta(R)}\exp\left(-n\inf_{z\in A\cap \mathcal{L}_n}\ent_{\utheta}(z)\right)\\
        & = (n+1)^{-d}\exp\left(-n\inf_{z\in A\cap \mathcal{L}_n}\rent_{\theta}(z)\right),
    \end{align*}
    and thereby also the counterparts of (2.1.14) and (2.1.15), i.e., 
    \begin{equation*}
        \limsup_{n\to\infty}\frac{1}{n}\log\sup_{\utheta\in B_\theta(R)}\probsub{\utheta}{\proc_n\in A} = -\liminf_{n\to\infty}\inf_{z\in A\cap \mathcal{L}_n}\rent_{\theta}(z)
    \end{equation*}
    and
    \begin{equation*}
        \liminf_{n\to\infty}\frac{1}{n}\log\sup_{\utheta\in B_\theta(R)}\probsub{\utheta}{\proc_n\in A}= -\limsup_{n\to\infty}\inf_{z\in A\cap \mathcal{L}_n}\rent_{\theta}(z).
    \end{equation*}
    Then, the upper bound of \eqref{eq_robust_Sanov} follows, as $A\cap\mathcal{L}_n\subseteq A$ for all $n$, so that  \[\liminf_{n\to\infty}\inf_{z\in A\cap \mathcal{L}_n}\rent_{\theta}(z)\geq\inf_{z\in A}\rent_{\theta}(z).\]

    Finally, to prove the lower bound of~\eqref{eq_robust_Sanov}, we fix an arbitrary $z \in \interior{A}$ and use the same sequence $\seq{z}{n} \subseteq \sspa$ as in the proof of~\cite[Theorem~2.1.10]{iB1}, which satisfies $z_n \in A \cap \mathcal{L}_n$ and $z_n \to z$ as $n \to \infty$. Then, since $\rent$ is continuous,
    \begin{equation*}
        -\limsup_{n\to\infty}\inf_{z\in A\cap \mathcal{L}_n}\rent_{\theta}(z)\geq-\limsup_{n\to\infty}\rent_{\theta}(z_n) = -\rent_{\theta}(z),
    \end{equation*}
    so that,  as $z$ was arbitrary from $\interior{A}$,
    \begin{equation*}
        \liminf_{n\to\infty}\frac{1}{n}\log\sup_{\utheta\in B_\theta(R)}\probsub{\utheta}{\proc_n\in A}\geq-\inf_{z\in \interior{A}}\rent_{\theta}(z).\qedhere
    \end{equation*}
\end{proof}

\bibliographystyle{agsm} 
\bibliography{bibfile}{}

\end{document}